\documentclass[11pt]{article}
\usepackage{amsmath,amssymb}
\usepackage[latin1]{inputenc}
\usepackage{epsfig}
\newtheorem{thm}{Theorem}

\newtheorem{lem}[thm]{Lemma}
\newtheorem{prop}[thm]{Proposition}

\newtheorem{rem}[thm]{Remark}
\newenvironment{dem}[1][]%
   {\ \\ {\bf Proof #1~: }}%
   {\hfill\mbox{\rule{2 true mm}{3 true mm}}\vskip 2 ex\noindent}
\newcommand{\E}{{\mathbb E}}
\renewcommand{\P}{{\mathbb P}}
\newcommand{\R}{{\mathbb R}}
\newcommand{\N}{{\mathbb N}}

\newcommand{\sm}{{s^-}}

\def\be{\begin{eqnarray}}
\def\ee{\end{eqnarray}}
\def\ben{\begin{eqnarray*}}
\def\een{\end{eqnarray*}}
\title{Nonlinear SDEs driven by L\'evy processes and related PDEs}
\author{Benjamin Jourdain\thanks{CERMICS, \'Ecole des Ponts, ParisTech, 6-8 avenue Blaise Pascal, Cit\'e
Descartes, Champs sur Marne, 77455 Marne la Vall\'ee Cedex 2,
e-mail:jourdain@cermics.enpc.fr}, Sylvie M\'el\'eard\thanks{CMAP,
Ecole Polytechnique, CNRS, route de Saclay, 91128 Palaiseau Cedex
e-mail: sylvie.meleard@polytechnique.edu}, Wojbor A.
Woyczynski\thanks{Department of Statistics and Center for
Stochastic and Chaotic Processes in Science and Technology, Case
Western Reserve University, Cleveland, OH 44106, e-mail:
waw@po.cwru.edu}}

\begin{document}

\maketitle

\begin{abstract}
In this paper we study general nonlinear
  stochastic differential equations, where the usual Brownian motion is replaced by  a L\'evy process. We also suppose that the
  coefficient multiplying the increments of this process is merely Lipschitz continuous and not necessarily linear
  in the time-marginals of the solution as is the case in the classical
  McKean-Vlasov model. We first study existence, uniqueness and particle
  approximations for these stochastic differential equations. When the
  driving process is a pure jump L\'evy process with a smooth but
  unbounded L\'evy measure, we develop a stochastic calculus of
  variations to prove that the time-marginals of the
  solutions are absolutely continuous with respect to the Lebesgue
  measure. In the case of a symmetric stable driving process, we deduce
  the existence of a function solution to a nonlinear
  integro-differential equation involving the fractional Laplacian.\par
 \par

 \end{abstract}


\noindent This paper studies  the following nonlinear stochastic differential
equation:
\begin{equation}
\begin{cases}
   X_t=X_0+\int_0^t\sigma(X_\sm,P_s)dZ_s,\quad t\in[0,T],\\
\forall s\in [0,T],\;P_s\;\mbox{ denotes the probability
  distribution of }\;X_s.
\end{cases}\label{eds}
\end{equation}
We assume that $X_0$ is a random variable with values in $\R^k$,
distributed according to $m$,
$(Z_t)_{t\leq T}$ a L\'evy process with values in $\R^d$,
independent of $X_0$, and $\sigma:\R^k\times{\mathcal P}(\R^k)\rightarrow
\R^{k\times d}$, where ${\mathcal P}(\R^k)$ denotes the set of
probability measures on $\R^k$. Notice that  the classical McKean-Vlasov model, studied for instance in
  \cite{Sznitman:91}, is obtained as a special case of \eqref{eds} by choosing
  $\sigma$ linear in the second variable and $Z_t=(t,B_t)$, with $\ B_t\ $ being a $(d-1)$-dimensional
  standard Brownian motion.

The first section of the paper is devoted to the existence problem and particle approximations for
\eqref{eds}.
Initially, we  address the case of square integrable both,  the  initial condition
$X_0$, and the L\'evy process $(Z_t)_{t\leq T}$.  Under these  assumptions the existence and uniqueness problem for
\eqref{eds} can be handled exactly as  in the Brownian case $Z_t=(t,B_t)$.  The nonlinear
stochastic differential equation \eqref{eds} admits a unique
solution as soon as $\sigma$ is Lipschitz continuous on
   $\R^k\times{\mathcal P}_2(\R^k)$  endowed with the product of the
   canonical metric on $\R^k$ and the Vaserstein metric $d$ on the set
   ${\mathcal P}_2(\R^k)$ of probability measures with finite second
   order moments. This assumption is much weaker than the assumptions imposed  on $\sigma$
   in the classical McKean-Vlasov model, where  it is also supposed
   to be linear in
   its second variable, that is,
   $\sigma(x,\nu)=\int_{\R^k}\varsigma(x,y)\nu(dy)$, for
   a Lipschitz continuous
   function  $\varsigma:\R^k\times\R^k\rightarrow \R^{k\times d}$. Then, replacing the  nonlinearity by
the related interaction, we define systems of $n$ interacting particles. In
the limit $n\rightarrow +\infty$, we prove, by a trajectorial
propagation of chaos result, that the dynamics of each particle
approximates the one given by \eqref{eds}. Unlike in the very
specific McKean-Vlasov model, where the universal $C/\sqrt{n}$ rate of convergence
 corresponds to the central limit theorem,
under our general assumptions on $\sigma$,  the
rate of convergence turns out to depend on the spatial dimension $k$.

In the next step,   the square
   integrability assumption is relaxed. However, to compensate for its loss,  we assume a reinforced
   Lischitz continuity of  $\sigma$ :  the
   Vaserstein metric $d$ on ${\mathcal P}(\R^k)$ is replaced by its smaller and
   bounded modification
   $d_1$ defined below. Then, choosing square integrable approximants of  the initial variable and the L\'evy process,  we prove existence for
   \eqref{eds}. Uniqueness remains an open question.

In the second section, we deal with the issue of absolute continuity of $P_t$ when
$Z$ is a pure jump L\'evy process with infinite intensity. For the sake of
simplicity of the exposition, we restrict ourselves to the one-dimensional case
$k=d=1$. When $\sigma$ does not vanish and admits two bounded derivatives
with respect to its first variable,  and the L\'evy measure of $Z$ satisfies some technical conditions, we prove that, for each $t>0$, $P_t$
has a density with respect to the Lebesgue measure on $\R$. The proof depends on
   a stochastic calculus of
variations for the SDEs driven by $Z$ which we develop by  generalizing the approach
of  Bichteler-Jacod \cite{Bichteler:83}, (see also Bismut
\cite{Bismut:83}), who dealt with the case of  homogeneous processes with a jump
measure equal to the Lebesgue measure. In our case,  the
nonlinearity induces an inhomogeneity in time and the jump
measure is much more general, which introduces additional
difficulties making the extension nontrivial. Graham-M\'el\'eard \cite{Graham:99}  developed similar techniques for a very specific stochastic
differential equation
  related to the Kac equation. In that case, the jumps
of the process were bounded. In our case, unbounded jumps are
allowed and we deal with the resulting possible lack of integrability of the
process $X$  by an appropriate conditioning.

In the third section, we keep the assumptions made on $\sigma$ in the second section, and
assume that the driving L\'evy process $Z$ is symmetric and $\alpha$-stable. Then, we  apply the absolute continuity results obtained in Section
2  to prove that  the solutions to \eqref{eds} are such that for
$t>0$, $P_t$ admits a density $p_t$ with respect to the Lebesgue measure on
the real line. In addition, calculating explicitly the adjoint of the generator of
$X$, we conclude that the function $p_t(x)$ is a weak solution to the
nonlinear Fokker-Planck equation
\begin{equation*}
   \begin{cases}
      \partial_tp_t(x)=D^\alpha_x(|\sigma(.,p_t)|^\alpha p_t(.))(x)\\
\lim_{t\rightarrow 0^+}p_t(x)dx=m(dx),
   \end{cases},
\end{equation*}
where, by a slight abuse of the notation, $\sigma(.,p_t)$ stands for
$\sigma(.,p_t(y)dy)$,  the limit is understood in the sense of the narrow convergence,
and $D^\alpha_x=-(-\Delta)^{\alpha/2}$ denotes the spatial, spherically symmetric fractional derivative of
order $\alpha$ defined here as a singular integral operator,
$$
D^\alpha_x f(x)=K\int_{\R} \left (f(x+y)-f(x)-{\bf 1}_{\{|y|\leq 1\}}f'(x)y
\right )\frac{dy}{|y|^{1+\alpha}},
$$
where $K$ is a positive constant.
For
$$\sigma(x,\nu)=\left(g_\varepsilon*\nu(x)\right)^s\mbox{
  with }\varepsilon>0,\;g_\varepsilon(x)=\frac{1}
  {\sqrt{2\pi\varepsilon}}e^{-\frac{x^2}{2\varepsilon}}\;\mbox{and}\;s>0,$$
one obtains the nonlocal approximation $\partial_t
p_t=D^\alpha_x((g_\varepsilon*p_t)^{\alpha s}p_t)$ of the fractional porous medium equation
$\partial_t p_t=D^\alpha_x(p_t^{\alpha s+1})$,
the physical interest of which is discussed at the end of the paper.
Other nonlinear evolution equations involving generators of L\'evy processes,
such as fractional conservation laws have been studied via probabilistic tools in, e.g., \cite{JMW:05a},
and\cite{JMW:05b}.

{\bf Notations :} Throughout the paper,  $C$ will denote a
constant which may change from line to line. In spaces with finite
dimension, the Euclidian norm is denoted by $|\;\;|$. Let
${\mathcal P}(\R^k)$ denote the set of probability measures on
$\R^k$, and ${\mathcal P}_2(\R^k)$ -- the subset of measures with
finite second order moments.  For $\mu,\nu\in{\mathcal
  P}_2(\R^k)$, the {\it Vaserstein metric}   is defined by the formula,
$$
d(\mu,\nu)=\inf\left\{\left(\int_{\R^k\times\R^k}|x-y|^2\,Q(dx,dy)\right)^{1/2}:Q\in{\mathcal
  P}(\R^k\times\R^k)\mbox{ with marginals }\mu\mbox{ and }\nu\right\}.
  $$
It induces the topology of weak convergence together with convergence of
  moments up to order $2$.
The modified Vaserstein metric on ${\mathcal
  P}(\R^k)$ defined by the formula,
$$
d_1(\mu,\nu)=\inf\left\{\left(\int_{\R^k\times\R^k}|x-y|^2\wedge 1\,Q(dx,dy)\right)^{1/2}:Q\in{\mathcal
  P}(\R^k\times\R^k)\mbox{ with marginals }\mu\mbox{ and }\nu\right\},
  $$
simply  induces the topology of weak convergence.


\section{Existence of a nonlinear process}
We first address the case when both, the initial condition $X_0$
and  $Z$ are square integrable, before relaxing these
integrability conditions later on.


\subsection{The square integrable case}
In this subsection we assume that the initial condition $X_0$, and
the L\'evy process $(Z_t)_{t\leq T}$, are both square integrable  :
$\E(|X_0|^2+|Z_T|^2)<+\infty$. Under this assumption, the
following inequality generalizes the Brownian case (see, \cite{Pr}, Theorem 66, p.339) :

\begin{lem}\label{doob}
   Let $p\geq 2$ be such that $\E(|Z_T|^p)<+\infty$. There is a constant $C_p$ such that, for any $\R^{k\times d}$-valued process
   $(H_t)_{t\leq T}$ predictable for the filtration ${\mathcal
   F}_t=\sigma(X_0,(Z_s)_{s\leq t})$, $\forall t\in[0,T]$,
$$
\E\left(\sup_{s\leq t}\left|\int_0^sH_udZ_u\right|^p\right)\leq C_p\int_0^t\E(|H_s|^p)ds.
$$
\end{lem}

Because of this inequality for $p=2$, the results obtained for the
classical McKean-Vlasov model driven by a standard Brownian motion
still hold. First, we state and prove

\begin{prop}
   Assume that $X_0$ and $(Z_t)_{t\leq T}$ are square integrable, and
   that the mapping $\sigma$ is Lipschitz continuous when
   $\R^k\times{\mathcal P}_2(\R^k)$ is endowed with the product of the
   canonical topology on $\R^k$ and the Vaserstein metric $d$ on
   ${\mathcal P}_2(\R^k)$. Then equation \eqref{eds} admits a unique
   solution such that $\E\left(\sup_{t\leq
       T}|X_t|^2\right)<+\infty$. Moreover, if for some $p>2$,
   $\E(|X_0|^p+|Z_T|^p)<+\infty$, then $\E\left(\sup_{t\leq
       T}|X_t|^p\right)<+\infty$.
\label{exuneds}\end{prop}

\begin{dem} We generalize here the pathwise fixed point approach well known
 in the classical  McKean-Vlasov case (see, Sznitman \cite{Sznitman:91}).
Let $\mathbb{D}$ denote the space
   of c\`adl\`ag
   functions from $[0,T]$ to $\R^k$, and ${\mathcal P}_2(\mathbb{D})$ the space of
   probability measures $Q$ on $\mathbb{D}$ such that $\int_\mathbb{D}\sup_{t\leq T}|Y_t|^2Q(dY)<+\infty$.
   Endowed with the Vaserstein metric
   $D_T(P,Q)$ where, for $t\leq T$, $$D_t(P,Q)=\inf\left\{\left(\int_{\mathbb{D}\times \mathbb{D}}
   \sup_{s\leq t}|Y_s-W_s|^2\,R(dY,dW)\right)^{1/2}:R\in{\mathcal
  P}(\mathbb{D}\times \mathbb{D})\mbox{ with marginals }P\mbox{ and }Q\right\},$$
${\mathcal P}_2(\mathbb{D})$ is a complete space.

For $Q\in{\mathcal P}(\mathbb{D})$ with time-marginals $(Q_t)_{t\in[0,T]}$,
in view of
   Lebesgue's Theorem, the distance
$$
d(Q_t,Q_s)\leq \int_{\mathbb{D}}|Y_t-Y_s|^2Q(dY)
$$
converges to $0$, as $s$
   decreases to $t$ (respectively,  $d(Q_{t^-},Q_s)\leq
   \int_{\mathbb{D}}|Y_{t^-}-Y_s|^2Q(dY)$ converges to 0, as $s$ increases to $t$; here $Q_{t^-}=Q\circ Y_{t^-}^{-1}$ is
   the weak limit of $Q_s$ as $s\rightarrow t^-$). Therefore, the mapping
   $t\in[0,T]\rightarrow Q_t$ is c\`adl\`ag when ${\mathcal P}_2(\R^k)$ is
   endowed with the metric $d$. As a consequence, for fixed $x\in\R^k$,
  the mapping  $t\in[0,T]\rightarrow \sigma(x,Q_t)$ is c\`adl\`ag. Hence, by a
   multidimensional version of Theorem
   6, p. 249, in \cite{Pr}, the
  standard stochastic differential equation
$$
X^Q_t=X_0+\int_0^t\sigma(X^Q_\sm,Q_s)dZ_s,\;t\in[0,T]
$$
admits a unique solution.

Let $\Phi$ denote the mapping on ${\mathcal P}_2(\mathbb{D})$ which associates the
law of $X^Q$ with $Q$. Let us check that $\Phi$ takes its values
in ${\mathcal P}_2(\mathbb{D})$. For $K>0$, we set $\tau_K=\inf\{s\leq T:|X^Q_s|\geq K\}$. By
  Lemma \ref{doob} and the Lipschitz property of $\sigma$, one has
\begin{align*}
\E\left(\sup_{s\leq t}|X^Q_{s\wedge\tau_K}|^2\right)&\leq
C\left(\E(|X_0|^2)+\int_0^t\E\left(1_{\{s\leq \tau_K\}}|\sigma(X^Q_s,Q_s)-\sigma(0,\delta_0)|^2
+|\sigma(0,\delta_0)|^2\right)ds\right)\\
&\leq C\left(\E(|X_0|^2)+\int_0^t\E\left(\sup_{r\leq
s}|X^Q_{r\wedge\tau_K}|^2\right)ds+t\int_\mathbb{D}\sup_{t\leq
T}|Y_t|^2Q(dY)+t|\sigma(0,\delta_0)|^2\right).
\end{align*}
By
Gronwall's Lemma, one deduces that
$$
\E\left(\sup_{s\leq T}|X^Q_{s\wedge\tau_K}|^2\right)\leq
C\left(\E(|X_0|^2)+|\sigma(0,\delta_0)|^2+\int_\mathbb{D}\sup_{t\leq
T}|Y_t|^2Q(dY)\right),
$$
 where the constant $C$ does not depend on
$K$. Letting $K$ tend to $+\infty$, one concludes by Fatou's Lemma
that
\begin{equation}
   \int_\mathbb{D}\sup_{s\leq T}|Y_s|^2d\Phi(Q)(Y)=\E\left(\sup_{s\leq
    T}|X^Q_s|^2\right)\leq C\left(\E(|X_0|^2)+|\sigma(0,\delta_0)|^2+\int_\mathbb{D}\sup_{t\leq T}|Y_t|^2Q(dY)\right).
\label{majol2}\end{equation}
Observe that a process
$(X_t)_{t\in [0,T]}$,  such that $\E\left(\sup_{t\leq
  T}|X_t|^2\right)<+\infty$, solves \eqref{eds} if
  and only if its law is a fixed-point of $\Phi$. So, to complete the proof
   of the Proposition, it suffices to check that $\Phi$ admits a unique fixed point.\par

By a formal computation, which can be made rigorous by a
localization
  procedure similar to the   one utilized above, for $P,Q\in{\mathcal P}_2(\mathbb{D})$ one has
\begin{align*}
   \E\left(\sup_{s\leq t}|X^P_{s}-X^Q_{s}|^2\right)&\leq
   C\int_0^t\E(|\sigma(X^P_\sm,P_s)-\sigma(X^Q_\sm,Q_s)|^2)ds\\&\leq
   C\int_0^t\E\left(\sup_{r\leq s}|X^P_{r}-X^Q_{r}|^2\right)+d^2(P_s,Q_s)ds.
   \end{align*}
By Gronwall's Lemma, one deduces that, $\forall t\leq T,$
$$
\;\E\left(\sup_{s\leq
    t}|X^P_{s}-X^Q_{s}|^2\right)\leq
    C\int_0^td^2(P_s,Q_s)ds.
    $$
Since $D^2_t(\Phi(P),\Phi(Q))\leq \E(\sup_{s\leq
  t}|X^P_s-X^Q_s|^2)$, and $d(P_s,Q_s)\leq D_s(P,Q)$, the last inequality implies,  $\forall t\leq T,$
$$
D^2_t(\Phi(P),\Phi(Q))\leq C\int_0^tD^2_s(P,Q)ds.
$$
Iterating this inequality, and denoting by $\Phi^N$ the $N$-fold
composition of $\Phi$, we obtain that, $\forall N\in\N^*, $
$$
D^2_T(\Phi^N(P),\Phi^N(Q))\leq
C^N\int_0^T\frac{(T-s)^{N-1}}{(N-1)!}D^2_s(P,Q)ds\leq
\frac{C^NT^N}{N!}D^2_T(P,Q).
$$
Hence, for $N$ large enough, $\Phi^N$
is a contraction which entails that
$\Phi$ admits a unique fixed point.

If, for some $p>2$,  $\E(|X_0|^p+|Z_T|^p)<+\infty$, a reasoning
similar to the one used in  the derivation of \eqref{majol2}, easily leads to the conclusion that
the constructed solution $(X_t)_{t\leq T}$ of \eqref{eds} is such
that
$$
\E\left(\sup_{s\leq
    T}|X_s|^p\right)\leq C\left(\E(|X_0|^p)+|\sigma(0,\delta_0)|^p+\E\left(\sup_{s\leq T}|X_s|^2\right)^{p/2}\right)<+\infty.
    $$
    \end{dem}

Our next step is  to study pathwise particle approximations for the
nonlinear process.
Let $((X^i_0,Z^i))_{i\in\N^*}$ denote a sequence of independent
pairs with $(X^i_0,Z^i)$ distributed like $(X_0,Z)$. For each
$i\geq 1$, let $(X^i_t)_{t\in[0,T]}$ denote the solution given by
Proposition \ref{exuneds} of  the nonlinear stochastic
differential equation starting from $X^i_0$ and driven by $Z^i$ :
\begin{equation}
\begin{cases}
   X^i_t=X^i_0+\int_0^t\sigma(X^i_\sm,P_s)dZ^i_s,\;t\in[0,T]\\
\forall s\in [0,T],\;P_s\;\mbox{ denotes the probability
  distribution of }\;X^i_s
\end{cases}.
\label{edsi}\end{equation} Replacing the nonlinearity by
interaction, we introduce the following system of $n$ interacting
particles
\begin{equation}
\begin{cases}
   X^{i,n}_t=X^i_0+\int_0^t\sigma(X^{i,n}_\sm,\mu^n_\sm)dZ^i_s,\quad t\in[0,T], \quad 1\le i\le n,\\
\mbox{ where }\mu^n=\frac{1}{n}\sum_{j=1}^n\delta_{X^{j,n}}\mbox{
  denotes the empirical measure}
\end{cases}.\label{systpart}
\end{equation}
Since for $\xi=(x_1,\hdots,x_n)$, and $\zeta=(y_1,\hdots,y_n)$ in
$\R^{nk}$, one has
\begin{equation}
 d\left(\frac{1}{n}\sum_{j=1}^n\delta_{x_j},\frac{1}{n}\sum_{j=1}^n\delta_{y_j}\right)\leq
\left(\frac{1}{n}\sum_{j=1}^n|x_j-y_j|^2\right)^{1/2}=\frac{1}{\sqrt{n}}|\xi-\zeta|.\label{vasdis}
\end{equation}
Existence of a unique solution to \eqref{systpart}, with finite second order moments,  follows from
Theorem 7, p. 253, in \cite{Pr}. Our next result establishes  the
trajectorial propagation of chaos result for the interacting particle system \eqref{systpart}.

\begin{thm}Under the assumptions of Proposition \ref{exuneds},
$$
\lim_{n\rightarrow +\infty}\sup_{i\leq n}\E\left(\sup_{t\leq T}|X^{i,n}_t-X^i_t|^2\right)=0
$$
Moreover, under additional assumptions, the following two  explicit estimates
hold :
\begin{itemize}
   \item If $\E(|X_0|^{k+5}+|Z_T|^{k+5})<+\infty$, then
\begin{equation}
   \sup_{i\leq n}\E\left(\sup_{t\leq
T}|X^{i,n}_t-X^i_t|^2\right)\leq Cn^{-\frac{2}{k+4}};\label{chaosmom}
\end{equation}
\item If
$\sigma(x,\nu)=\int_{\R^k}\varsigma(x,y)\nu(dy)$,  where
  $\varsigma:\R^k\times\R^k\rightarrow\R^{k\times d}$ is a Lipschitz
  continuous function, then
\begin{equation}
   \sup_{i\leq n}\E\left(\sup_{t\leq T}|X^{i,n}_t-X^i_t|^2\right)\leq
\frac{C}{n},\label{chaosmckean}\end{equation}
\end{itemize}

where the constant $C$ does not depend on
$n$.\label{chaos}
\end{thm}

\noindent The proof of the first assertion relies
on the following

\begin{lem}\label{majovasind}
 Let $\nu$ be
a probability measure on $\R^k$ such that
$\int_{\R^k}|x|^2\nu(dx)<+\infty\ $, and
$\nu^n=\frac{1}{n}\sum_{j=1}^n\delta_{\xi_j}$ denote the empirical
measure associated with a sequence $(\xi_i)_{i\geq 1}$ of
independent random variables with law $\nu$. Then, $\forall n \geq 1,$
$$
\E\left(d^2\left(\nu^n,\nu\right)\right)\leq
4\int_{\R^k}|x|^2\nu(dx),\;\;\mbox{and}\;\;\lim_{n\rightarrow
  +\infty}\E\left(d^2\left(\nu^n,\nu\right)\right)=0.
  $$
\end{lem}

\begin{dem}[of Lemma \ref{majovasind}]
By the strong law of large numbers,  as $n$ tends to $\infty$, almost surely $\nu^n$
converges weakly to $\nu$ and, $\forall i,j\in\{1,\hdots,k\}$,
$\int_{\R^k}x_i\nu_n(dx)$ (resp. $\int_{\R^k}x_ix_j\nu_n(dx)$)
converges to $\int_{\R^k}x_i\nu(dx)$ (resp.
$\int_{\R^k}x_ix_j\nu(dx)$). Since the
Vaserstein distance $d$ induces the topology of simultaneous weak
convergence and convergence of moments up to order $2$, one
deduces that almost surely, $d(\nu^n,\nu)$ converges to $0$, as $n$
tends to $\infty$. Hence, to
conclude the
proof of the first assertion, it is enough to check that the random
variables $(d^2(\nu^n,\nu))_{n\geq 1}$ are uniformly integrable. To see that  note the inequality
$$
d^2(\nu^n,\nu)\leq\frac{2}{n}\sum_{j=1}^n
|\xi_j|^2+2\int_{\R^k}|x|^2\nu(dx).
$$
 The right-hand side is
nonnegative and converges almost surely to
$4\int_{\R^k}|x|^2\nu(dx)$, as $n\rightarrow \infty$. Since its
expectation is constant, and equal to the expectation  of the limit, one
deduces that the convergence is also  in $L^1$. As a consequence, for $n\geq
1$, the random variables in the right-hand side, and therefore in
the left-hand side, are
uniformly integrable.
\end{dem}

\eject

\begin{dem}[of Theorem \ref{chaos}]
Let $P^n=\frac{1}{n}\sum_{j=1}^n\delta_{X^j}$
  denote the empirical measure of the independent nonlinear
  processes \eqref{edsi}. By a formal computation, which can be made rigorous by a localization
argument similar to the one made in the proof of Proposition
\ref{exuneds}, one has, $  \forall t\leq T,$
\begin{align*}
 \E\left(\sup_{s\leq t}|X^{i,n}_s-X^i_s|^2\right)\leq
    &C\int_0^t\E\left(|\sigma(X^{i,n}_s,\mu^n_s)-\sigma(X^i_s,P^n_s)|^2\right)ds\\
    &+C\int_0^t\E\left(|\sigma(X^i_s,P^n_s)-\sigma(X^{i}_s,P_s)|^2\right)ds.
\end{align*}
In view of  the Lipschitz property of $\sigma$, the estimate \eqref{vasdis}, and the
exchangeability of the couples $(X^i,X^{i,n})_{1\leq i\leq n}$,
the first term of the right  is smaller than
$
C\int_0^t\E\left(\sup_{r\leq s}|X^{i,n}_r-X^i_r|^2\right)ds.
$
By
Gronwall's Lemma, and the Lipschitz assumption on $\sigma$, one
deduces that
$$
\E\left(\sup_{t\leq T}|X^{i,n}_t-X^i_t|^2\right)\leq
   C\int_0^T\E\left(\left|\sigma\left(X^i_s,P^n_s\right)-\sigma(X^i_s,P_s)\right|^2\right)ds\leq
 C\int_0^T\E(d^2(P^n_s,P_s))ds.
 $$
The first assertion then follows from Lemma \ref{majovasind}, the
upper-bounds of the second order moments given in Proposition
\ref{exuneds}, and by Lebesgue's Theorem.

The second assertion is
deduced from the upper-bounds for  moments of order $k+5$ combined
with the following restatement of Theorem 10.2.6 in \cite{raru} :
$$
\E\left(d^2\left(P_s^n,P_s\right)\right)\leq
C\left(1+\sqrt{\int_{\R^k}|y|^{k+5}P_s(dy)}\right)n^{-\frac{2}{k+4}},
$$
where the constant $C$ only depends on $k$. The precise dependence of the upper-bound on
$\int_{\R^k}|y|^{k+5}P_s(dy)$ comes from a carefull reading of
the proof given in \cite{raru}.

Finally, if, as in the usual McKean-Vlasov framework (see
\cite{Sznitman:91}),
$\sigma(x,\nu)=\int_{\R^k}\varsigma(x,y)\nu(dy)$, where
  $\varsigma=(\varsigma_{ab})_{a\leq k,b\leq d}:\R^k\times\R^k\rightarrow\R^{k\times d}$ is Lipschitz
  continuous, then  $\E\left(\left|\sigma\left(X^i_s,P^n_s\right)-\sigma(X^i_s,P_s)\right|^2\right)$
  is equal to
\begin{align*}
  \sum_{a=1}^k\sum_{b=1}^d\frac{1}{n^2}\sum_{j,l=1}^n\E\left(\left[\varsigma_{ab}(X^i_s,X^j_s)
  -\int_{\R^k}\varsigma_{ab}(X^i_s,y)P_s(dy)\right]\left[\varsigma_{ab}(X^i_s,X^l_s)-\int_{\R^k}\varsigma_{ab}(X^i_s,y)P_s(dy)\right]\right).
\end{align*}
Since, by independence of the random variables $X^1_s,\hdots,X^n_s$
with common law $P_s$, the expectation in the above summation
vanishes as soon as $j\neq l$, the third assertion of Theorem 3 easily follows.
\end{dem}

\begin{rem}$ $ \rm
\begin{itemize}
\item \label{buckwise} \rm Observe  the lower estimate
$$
d(\nu^n,\nu)\geq \left(\int_{\R^k}\min_{1\leq j\leq
  n}|\xi_j-x|^2\nu(dx)\right)^{1/2}\geq \inf_{(y_1,\hdots,y_n)\in(\R^{k})^n}\left(\int_{\R^k}\min_{1\leq j\leq
  n}|y_j-x|^2\nu(dx)\right)^{1/2}.
  $$
Moreover, according to the Bucklew and Wise Theorem
\cite{buckwise}, if $\nu$ has a density $\varphi$
  with respect to the Lebesgue measure on $\R^k$ which belongs to $L^{\frac{k}{2+k}}(\R^k)$, then,
as $n$ tends to infinity,
$$
n^{1/k}\inf_{(y_1,\hdots,y_n)\in(\R^{k})^n}\left(\int_{\R^k}\min_{1\leq
j\leq
  n}|y_j-x|^2\varphi(x)dx\right)^{1/2}
  $$
  converges to
  $C_k\|\varphi\|_{\frac{k}{2+k}}$, where the constant $C_k$ only depends
  on $k$. Hence, one cannot expect $\E(d^2(\nu^n,\nu))$ to vanish
  quicker than $Cn^{-2/k}$.

Therefore, if  $\nu\rightarrow\sigma(x,\nu)$ is merely Lipschitz
continuous
  for the Vaserstein metric, one cannot expect $\E\left(\sup_{t\leq
      T}|X^{i,n}_t-X^i_t|^2\right)$ to vanish quicker than
  $Cn^{-2/k}$. The rate of convergence obtained in \eqref{chaosmom} is
  not far from being optimal at least for a large spatial dimension $k$.
Nevertheless, in the
McKean-Vlasov framework, where the structure of $\sigma$ is very specific, one can overcome this dependence of the convergence rate   on the dimension $k$, and recover the usual central limit
theorem rate.

\item The square integrability assumption on the initial variable $X_0$
can be relaxed if $\sigma$ is Lipschitz continuous with
$\R^k\times{\mathcal P}(\R^k)$  endowed with the product of the
canonical topology on $\R^k$ and the modified Vaserstein metric
$d_1$ on ${\mathcal P}(\R^k)$. Indeed, one may then adapt the
fixed-point approach in the proof of Proposition \ref{exuneds} by
defining ${\mathcal P}$ as the space of probability measures on
$\mathbb{D}$, and replacing $\sup_{s\leq t}|Y_s-W_s|^2$ by
$\sup_{s\leq t}|Y_s-W_s|^2\wedge 1$ in the definition of
$D_t(P,Q)$. This way, one obtains that the nonlinear stochastic
differential equation \eqref{eds} still admits a unique solution
even if the initial condition $X_0$ is not square integrable.
Moreover, for any probability measure $\nu$ on $\R^k$, if $\nu_n$
is defined as above, $\E(d_1^2(\nu_n,\nu))$ remains bounded by one,
and converges to $0$ as $n$ tends to infinity. Therefore, the first
assertion in Theorem \ref{chaos} still holds.

\end{itemize}

\end{rem}

\noindent The next subsection  is devoted to the more complicated case when
the square integrability assumption on the L\'evy process
$(Z_t)_{t\leq T}$ is also relaxed.


\subsection{The general case}
In this section, we impose no integrability conditions,
either on the initial condition $X_0$, or on the L\'evy process
$(Z_t)_{t\leq
  T}$. Let us denote by $m\in{\mathcal P}(\R^k)$ the distribution of the
former. According to the L\'evy-Khintchine formula, the
infinitesimal generator of the latter can be written,  for $f\in C^2_b(\R^d)$, in the form
$$
Lf(z)=\frac{1}{2}\sum_{i,j=1}^d
a_{ij}\partial^2_{z_i,z_j}f(z)+b.\nabla
f(z)+\int_{\R^d}\left[f(z+y)-f(z)-{\bf
  1}_{\{|y|\leq 1\}}y.\nabla f(z)\right]\beta(dy),
  $$
where $a=(a_{ij})_{1\leq i,j\leq d}$ is a non-negative symmetric
matrix, $b$ a given vector in $\R^d$, and $\beta$ a measure on
$\R^d$ satisfying the integrability condition
$\int_{\R^d}(1\wedge|y|^2)\beta(dy)<+\infty$.

To deal with non square integrable sources of randomness, we
impose  a stronger continuity condition on $\sigma$, namely,  we assume
that $\sigma$ is Lipschitz continuous when $\R^k\times{\mathcal
P}(\R^k)$ is endowed with the product of the canonical topology on
$\R^k$ and the modified Vaserstein metric $d_1$ on ${\mathcal
P}(\R^k)$.  Notice that, under this assumption, for each
$x\in\R^k$, the mapping $\nu\in{\mathcal P}(\R^k)\rightarrow
\sigma(x,\nu)$ is bounded.

In order to prove existence of a weak
solution to \eqref{eds}, we introduce a cutoff parameter
$N\in{\mathbb N}^*$, and define a square integrable initial random
variable $X_0^N=X_01_{\{|X_0|\leq N\}}$, and a square integrable
L\'evy process $(Z^N_t)_{t\leq T}$ by removing the jumps of
$(Z_t)_{t\leq T}$ larger than $N$ :
$$
Z^N_t=Z_t-\sum_{s\leq
t}1_{\{|\Delta
  Z_s|> N\}}\Delta Z_s.
  $$
Let $(X^N_t)_{t\in[0,T]}$ denote the solution given by Proposition
\ref{exuneds} of the nonlinear stochastic differential equation
starting from $X^N_0$ and driven by $(Z^N_t)_{t\in[0,T]}$ :
\begin{equation}
\begin{cases}
   X^N_t=X^N_0+\int_0^t\sigma(X^N_\sm,P^N_s)dZ^N_s,\;t\in[0,T]\\
\forall s\in [0,T],\;P^N_s\;\mbox{ denotes the probability
  distribution of }\;X^N_s
\end{cases}.\label{edsb}
\end{equation}
We are going to prove that when the cutoff parameter $N$ tends to
$\infty$, then $(X^N_t)_{t\in[0,T]}$ converges in law to a weak
solution of \eqref{eds}. More precisely, let us denote by $P^N$
the distribution of $(X^N_t)_{t\in[0,T]}$, and by $(Y_t)_{t\in
[0,T]}$ the canonical process on ${\mathbb D}$.

\begin{prop}\label{proppdmgene}
   The set of probability measures $(P^N)_{N\in{\mathbb N}^*}$ is tight when
   ${\mathbb D}$ is endowed with the Skorohod topology. In addition, any weak
   limit $P$, with time marginals $(P_t)_{t\in [0,T]}$, of its  converging
   subsequences solves the following martingale problem :
\begin{equation}
   \begin{cases}
      P_0=m\;\;\mbox{and}\;\;\forall \varphi:\R^k\rightarrow
      \R,\;C^2\;\mbox{with compact support},\\
\left(M^\varphi_t=\varphi(Y_t)-\varphi(Y_0)-\int_0^t{\mathcal
     L}[P_s]\varphi(Y_s)ds\right)_{t\in[0,T]}\;\mbox{is a $P$-martingale}
   \end{cases}\label{pdm},
\end{equation}
where for each $\nu\in{\mathcal P}(\R^k)$, and any $x \in \R^k$,
\begin{align}
   {\mathcal
     L}[\nu]\varphi(x)=&\frac{1}{2}\sum_{i,j=1}^k(\sigma
a\sigma^*(x,\nu))_{ij}\partial^2_{x_i,x_j}\varphi(x)+(\sigma(x,\nu)b).\nabla\varphi(x)\notag\\
&+\int_{\R^d}\left[\varphi\left(x+\sigma(x,\nu)y\right)-\varphi(x)-{\bf
  1}_{\{|y|\leq 1\}}\sigma(x,\nu)y.\nabla\varphi(x)\right]\beta(dy)ds.
\label{gene}\end{align}
\end{prop}

\begin{dem}
   Let us first remark that for $N\in{\mathbb N}^*$, and  for a  fixed
   $x\in\R^k$, the mapping $t\in [0,T]\rightarrow \sigma(x,P^N_t)$ is
   c\`adl\`ag and bounded by a constant not depending on $N$. As a
   consequence, according to Theorem 6 p. 249 \cite{Pr}, for a  fixed
   $M\in\N^*$, the stochastic
   differential equation
$$
X^{N,M}_t=X^{N\wedge
  M}_0+\int_0^t\sigma(X^{N,M}_{s-},P^N_s)dZ^{N\wedge M}_s,\;\quad t\in[0,T],
  $$
admits a unique solution. Let us denote by $P^{N,M}$ the law of
  $(X^{N,M}_t)_{t\in[0,T]}$. By trajectorial uniqueness, $\forall N\in{\mathbb N}^*,\;\forall t\in [0,T],$
$$
X^N_t=X^{N,M}_t,
$$
as long as $|X_0|\vee\sup_{t\in]0,T]}|\Delta Z_t|\leq M$. The
probability of the latter event tends to one as $M$ tends to
infinity. Using both the necessary and the sufficient conditions
of Prokhorov's Theorem, one deduces that the tightness of the
sequence $(P^N)_{N\in\N^*}$ is implied by the tightness of the
sequence $(P^{N,M})_{N\in\N^*}$, for any fixed $M\in\N^*$.

\vskip -2mm

Let us
now prove this last
result by fixing $M\in\N^*$.
Using the boundedness of $d_1$, one
easily checks that
\begin{equation}
   \sup_{N\in\N^*}\E\left(\sup_{t\leq T}|X^{N,M}_t|^2\right)<+\infty.\label{majocarunifN}
\end{equation}
This implies tightness of the laws of the random variables
$\left(\sup_{t\leq
  T}|X^{N,M}_t|\right)_{N\in\N^*}$. In order to use Aldous' criterion, we set
  $\varepsilon,\delta>0$, and
  introduce two stopping times $S$, and $\tilde{S}$, such
  that $0\leq S\leq \tilde{S}\leq (S+\delta)\wedge T$. Let us also
  remark that, for $K\in{\mathbb N}^*$, and $b^K=b+\int_{\R^d}y{\bf
  1}_{\{1<|y|\leq K\}}\beta(dy)$, the process $\left(\tilde{Z}^K_t=Z^K_t-b^Kt\right)_{t\in[0,T]}$
  is a centered L\'evy process and therefore a martingale. Now, observe that
\begin{align}
   \P\left(|X^{N,M}_{\tilde{S}}-X^{N,M}_S|^2\geq
   \varepsilon\right)\leq&\P\left(\left|\int_S^{\tilde{S}}\sigma(X^{N,M}_s,P^N_s)b^{N\wedge
   M}ds\right|^2\geq
   \frac{\varepsilon}{4}\right)\notag\\&+\P\left(\left|\int_S^{\tilde{S}}\sigma(X^{N,M}_{s^-},P^N_s)d\tilde{Z}^{N\wedge
   M}_s\right|^2\geq \frac{\varepsilon}{4}\right).\label{majoaccroiss}
\end{align}
Using the boundedness of the sequence $(b^{N\wedge
M})_{N\in\N^*}$, the Lipschitz property of $\sigma$ with respect
to its first variable,  and \eqref{majocarunifN} combined with the
inequalities of Markov and Cauchy-Schwarz, one obtains that the
first term of the right-hand-side is smaller than
$C\delta^2/\varepsilon$, where the constant $C$ does not depend on
$N$. For the second term of the right-hand-side, one remarks that
Doob's optional sampling Theorem, followed by  the Lipschitz property of
$\sigma$, and \eqref{majocarunifN}, imply that
\begin{align*}\E\left(\left|\int_S^{\tilde{S}}\sigma(X^{N,M}_{s^-},P^N_s)d\tilde{Z}^{N\wedge
   M}_s\right|^2\right)&=\E\bigg(\int_S^{\tilde{S}}\bigg[\sum_{i=1}^k(\sigma
   a\sigma^*)_{ii}(X^{N,M}_s,P^N_s)\\&\phantom{\E\bigg(\int_S^{\tilde{S}}\bigg[}+\int_{\R^d}|\sigma(X^{N,M}_s,P^N_s)y|^2{\bf
   1}_{\{|y|\leq N\wedge M\}}\beta(dy)\bigg]ds\bigg)\\
&\leq C\delta,
\end{align*}
where $C$ does not depend on $N$. By
Markov's Inequality, the second term
   of the right-hand-side of \eqref{majoaccroiss} is smaller than
   $C\delta/\varepsilon$ and, in view of Aldous' criterion,  we conclude that the sequence
   $(P^{N,M})_{N\in\N^*}$ is tight.

\vskip -2mm
Now, let us  denote by $P$ the limit of a converging subsequence of
   $(P^N)_{N\in\N^*}$ that we still index by $N$ for simplicity's sake. Also, let  $\varphi$ denote a compactly
   supported $C^2$ function on $\R^k$. For $p\in\N^*$, $0\leq s_1\leq s_2\leq \hdots\leq s_p\leq
   s\leq t\leq T$,  and a continuous and
   bounded function $\psi:(\R^k)^p\rightarrow \R$ , let $F$ denote the mapping on ${\mathcal P}({\mathbb
     D})$ defined by
$$
F(Q)=\int_{\mathbb
     D}\left(\varphi(Y_t)-\varphi(Y_s)-\int_s^t{\mathcal
     L}[Q_u]\varphi(Y_u)du\right)\psi(Y_{s_1},\hdots,Y_{s_p})Q(dY).
     $$
For $F^N$ defined like $F$, but  with ${\bf 1}_{\{|y|\leq N\}}\beta(dy)$
replacing
     $\beta(dy)$ in the definition \eqref{gene} of ${\mathcal L}[\nu]$, one has
     $F^N(P^N)=0$. Therefore
$$
|F(P^N)|=|F(P^N)-F^N(P^N)|\leq
2(t-s)\|\psi\|_{\infty}\|\varphi\|_{\infty}\int_{\R^d}{\bf
1}_{\{|y|\geq
  N\}}\beta(dy)\stackrel{N\rightarrow +\infty}{\longrightarrow}0.
  $$
The mapping $(x,\nu)\in \R^k\times{\mathcal P}(\R^k)\rightarrow
{\mathcal
     L}[\nu]\varphi(x)$ is bounded, continuous in $x$ for a  fixed $\nu$,
   and continuous in $\nu$, uniformly for $x\in\R^k$. Therefore, as soon as $s_1,\hdots,s_p,s,t$ do not belong to the at most
   countable set $\{u\in]0,T]:\;P(\Delta Y_u\neq 0)>0\}$, then $F$ is
   continuous and bounded at point $P$ which implies
   $F(P)=\lim_{N\rightarrow +\infty}F(P^N)=0$.   In view of  the right continuity of $u\rightarrow Y_u$ and Lebesgue's
   theorem, this
   equality still holds without any restriction on $s_1,\hdots,s_p,s,t$.
   Hence, $(M_t^\varphi)_{t\in[0,T]}$ is a $P$-martingale. Since
 the sequence $(X^N_0)_{N\in\N^*}$ converges in distribution to $X_0$,
 $P_0=m$,  which concludes the proof.\end{dem}

The above existence result for the
martingale problem   \eqref{pdm} implies an analogous existence statement  for the corresponding
nonlinear Fokker-Planck equation.

\begin{prop}\label{propedsgene}
Let $P$ denote a solution of \eqref{pdm}. Then the time marginals $(P_t)_{t \in[0,T]}$ solve the initial value problem
\begin{equation}
   \partial_t P_t={\mathcal L}^*[P_t]P_t,\quad P_0=m, \label{fokplan}
\end{equation}
in the weak sense,  where,  for $\nu\in {\mathcal P}(\R^k)$,
${\mathcal
  L}^*[\nu]$ denotes the formal adjoint of ${\mathcal
  L}[\nu]$ defined by the following condition:   $\forall \phi,\psi$ $C^2$ with
compact support on $\R^k$,
$$
\int_{\R^k}{\mathcal L}^*[\nu]\psi(x)\varphi(x)dx=\int_{\R^k}
\psi(x){\mathcal L}[\nu]\varphi(x)dx.
$$
 Moreover, the standard
stochastic differential equation
\begin{equation}
   X^P_t=X_0+\int_0^t\sigma(X^P_{s^-},P_s)dZ_s\label{linsde}
\end{equation}
admits a unique solution $(X^P_t)_{t\in[0,T]}$ and, for each
$t\in[0,T]$, $X^P_t$ is distributed according to $P_t$.
\end{prop}

\begin{dem}
   The first assertion follows readily from the constancy of the expectation of the
$P$-martingale $(M^\varphi_t)_{t\in[0,T]}$.
Existence and uniqueness for \eqref{linsde} follows from \cite{Pr},
Theorem 6 p. 249. Now, if $Q_t$ denotes the law of $X^P_t$ for
$t\geq 0$, then $(Q_t)_{t\geq 0}$ solves
$$
\partial_t Q_t={\mathcal L}^*[P_t]Q_t,\quad Q_0=m,
$$
 in the weak sense. Since
$(P_t)_{t\geq 0}$ also is a weak solution of   this linear equation, by
Theorem 5.2 \cite{bk}, one concludes that $\forall t\leq T$,
$P_t=Q_t$.
\end{dem}

\begin{rem} \rm
   We have not been able to prove uniqueness for the nonlinear martingale problem   \eqref{pdm}.
    However, our assumptions, and Theorem
   6,  p.249, in  \cite{Pr}, ensure existence and uniqueness for the particle
   system \eqref{systpart}. Like in the proof of Proposition \ref{proppdmgene}, one
   can check that the laws of the processes $X^{1,n},\;n\geq 1$, are
   tight. According to \cite{Sznitman:91}, this implies uniform tightness of the
   laws $\pi_n$ of the empirical measures $\mu^n$. For a  fixed $x\in\R^k$, the
   function $\nu\rightarrow\sigma(x,\nu)$ is continuous and bounded when
   ${\mathcal P}(\R^k)$ is endowed with the weak convergence
   topology. Then one can prove that the limit points of the sequence
   $(\pi_n)_n$ give full weight to the solutions of the nonlinear martingale
   problem \eqref{pdm}.
\end{rem}


\section{Absolute continuity of the marginals}
    In  this section we restrict ourselves to the one-dimensional case
    $k=d=1$, and assume that $Z$ is a pure jump L\'evy
    process with a L\'evy measure $\beta$ which admits a density, say $\beta_1$,  in  the neighborhood of the origin, that is
$$
\beta(dy)={\bf 1}_{|y|\leq
      1}\beta_1(y)dy+{\bf 1}_{|y|>
    1}\beta(dy).
    $$
We set  $\beta_1(y)=0$,  for $|y|>1$.
Then \begin{align}
    \label{Levy-dens}
    Z_t=\int_{(0,t]\times \R}y\;\tilde{N}_1(ds,dy)+ \int_{(0,t]\times \R}y\;
    N_2(ds,dy),
    \end{align}
    where $N_1$, and $N_2$, are two independent  Poisson point measures on
    $\R_+\times \R$ with intensity measures equal,  respectively,  to $1_{|y|\leq
      1}\beta_1(y)dy$, and $1_{|y|>
    1}\beta(dy)$, and $\tilde{N}_1$ is the compensated
    martingale measure of $N_1$.

\noindent We work here \begin{itemize}
   \item either under the assumptions of Subsection 1.1, i.e., ${\mathbb
    E}(|X_0|^2+|Z_T|^2)<+\infty$ and Lipschitz continuity of
$\sigma$ on $\mathbb{R}\times{\cal P}_2(\mathbb{R})$ endowed with
the product of the Euclidean metric and the Vaserstein
metric,
\item or with the general assumptions of Subsection 1.2, i.e., no integrability
conditions on $Z$ and  $X_0$, and Lipschitz continuity of
$\sigma$ on $\mathbb{R}\times{\cal P}(\mathbb{R})$ endowed with
the product of the Euclidean metric and the modified Vaserstein
metric $d_1$.
\end{itemize}

\noindent Propositions \ref{exuneds} and \ref{proppdmgene} ensure
the existence of a probability measure  solution $P$ of
\eqref{pdm}. According to Proposition \ref{propedsgene}, there is
a unique pathwise solution $X$ (which is then unique in law),
  to the   stochastic differential equation
\begin{equation}
\label{eds-stable}
  X_t=X_0+\int_{(0,t]\times \mathbb{R}}\sigma(X_\sm,P_s)\ y\tilde{N}_1(dy,ds)
   +\int_{(0,t]\times \mathbb{R}}\sigma(X_\sm,P_s)\ yN_2(dy,ds),
\end{equation}
and, for $t\in [0,T]$, $X_t$ is distributed according to the
time marginal $P_t$.

Roughly speaking, our goal  is to prove that, for each  $t>0$, the
    probability measure $P_t$ has a density with respect to the Lebesgue measure  as long as
 the measure $\beta$, restricted to
$[-1,1]$, has an infinite total
    mass  due to an  explosion of the density function $\beta_1(y)$ at $0$.
     Indeed, we have in this case an
accumulation of small jumps immediately after time $0$, which will
imply the absolute continuity of the law of $X_t$ under suitable
regularity assumptions on $\beta_1$.

 For this purpose we develop a stochastic calculus of
variations for diffusions with jumps driven by the L\'evy process
defined in  \eqref{Levy-dens}. We thus generalize the approach
developed by  Bichteler and Jacod \cite{Bichteler:83}  (also, see
Bismut \cite{Bismut:83}),  who considered homogeneous processes
with a jump measure equal to the Lebesgue measure. Here, the
nonlinearity introduces an inhomogeneity in time, and the jump
measure is much more general, which complicates the situation
considerably and introduces  additional difficulties.  Graham and
M\'el\'eard \cite{Graham:99}  developed similar techniques for a
very specific stochastic differential equation
  related to the Kac equation. In that case, the jumps
of the process were bounded.   In our case, unbounded jumps are
also allowed.

Our approach requires  that we make the following standing assumptions on
 the coefficient $\sigma(x,\nu)$, and the L\'evy density $\beta_1$:

\noindent \textbf{Hypotheses (H):}

{\bf 1.} The coefficient $\sigma(x,\nu)$ is twice differentiable in $x$.

{\bf 2.}  There exist constants  $K_1$, and $K_2$, such that,  for each $x,$ and  $\nu$,
\be \label{hypsigma} \ |\sigma'_x(x,\nu)| \leq K_1,  \quad  {\rm and } \quad
\quad|\sigma''_x(x,\nu)| \leq K_2. \ee

{\bf 3.}  For each $x,$ and $ \nu$,
 \be \sigma(x,\nu)\neq 0.\ee

\noindent \textbf{Hypotheses ($\bf H_1$):}

{\bf 1.}  The function $\beta_1$ is twice continuously differentiable
away from $\{0\}$.

{\bf 2.}  \be \label{explose} \int_{-1}^1 \beta_1(y)dy=+\infty. \ee

{\bf 3.} There exists a non-negative and non-constant function $k$ of class $C^1$ on
$[-1,1]$ such that $k(-1)=k(1)=0$, and such that
\begin{itemize}
\item \be \label{hypint} \int_{-1}^1 k^2(y) \beta_1(y)dy
<+\infty, \ee \item  for all $y\in [-1, 1]$,
 \be
 \label{hypborne}
  \sup_{a\in[-K_1,
K_1],\lambda\in[-1,1]}{1\over |\lambda|}\left|{\beta_1(y+\lambda(1+a
y)k(y))\over \beta_1(y)}(1+\lambda(a k(y)+(1+ay)k'(y))-1\right|
\leq {1\over 2},\nonumber\\
\ \ee \item
 \be
 \label{hypbeta}
  \sup_{a\in[-K_1,
K_1]}\int_{-1}^1\left|\frac{\beta_1'(y)}{\beta_1(y)}(1+ay)k(y)+ak(y)+(1+a y)k'(y)\right|^2 \beta_1(y)dy <
+\infty,\nonumber\\
\ \ee \item
 \be
 \label{hypderive}
  &&\sup_{a\in[-K_1,
K_1]}\int_{-1}^1\sup_{\lambda\in[-1,1]}\bigg(\left|{\beta_1''(y+\lambda(1+a y)k(y))\over \beta_1(y)}
(1+\lambda(a k(y)+(1+ay)k'(y))\right|^2k^2(y) \nonumber\\
&&\hskip 4cm +\left|{\beta_1'(y+\lambda(1+a
y)k(y))\over \beta_1(y)}\right|^2 \bigg)k^2(y) \beta_1(y)dy <
+\infty,\nonumber\\
\ \ee 
 \item
  for all $y\in [-1, 1]$,
  \be
  \label{hypbound} |k(y)| < {1\over 4(1+ K_1)}\quad ; \quad |k'(y)| <{1\over
  4(1+K_1)}.
  \ee
\end{itemize}

The assumption ($\bf H_13$)  on $\beta_1$ is obviously technical, but the assumption  ($\bf H_12$)  is essential, and cannot be avoided if one hopes to prove the absolute continuity result.
The main example  satisfying  assumptions ($\bf H_1$) is
    the symmetric stable process with index $\alpha\in(0,2)$, for
    which $\beta(dy)={K\,dy/ |y|^{1+\alpha}}$, as developed in Section 3.

\begin{thm}
 \label{density} Consider the real-valued  process $X$ satisfying
    the nonlinear stochastic differential equation
    \eqref{eds-stable}.
    Assume that   $\sigma$ satisfies Hypotheses {\rm  ($\bf H$)},  and that
     $\beta_1$ satisfies Hypotheses {\rm ($\bf H_1$)}.
     Then the law of the real-valued random  variable $X_T$
    has a density with respect to the Lebesgue measure.
    \end{thm}


The remainder  of this section is devoted to the proof of Theorem
\ref{density} which  will proceed through a series of lemmas and
propositions. Our aim is to show that $P_T$ has a density with
respect to the Lebesgue's measure. Because of the compensated
martingale term $\tilde{N}_1$  it would be natural to work with
square integrable processes. But  the finiteness of
$\E\left(\sup_{t\leq T}|X_t|^2\right)$ is not guaranteed  because
of   the big jumps of $N_2$ . So,  to develop the relevant
stochastic calculus of variations in $L^2$, we use a trick
defining $\P_T$ as the conditional law of $(X_0,N_1,N_2)$ given
$(X_0,N_2^T)$,  where $N_2^T$ denotes the restriction of the
measure $N_2$ to $[0,T]\times\R$. Thus, in what follows, the
random variables considered
 being functions of $(X_0,N_1,N_2)$ we may define their
$\E_T$-expectations as the integral of the corresponding functions
under $\P_T$. From now on, for notational simplicity, every
statement concerning $\E_T$, or $\P_T$, holds
almost everywhere under the law of $(X_0,N_2^T)$, even if this fact  is not mentioned explicitly.
This conditioning allows us to use the same
techniques as if the process $X$ were square integrable. More
precisely, given $N_2^T$, there are finitely many  jump times and jump amplitudes of
$N_2$ on $(0,T]$ and we will denote them  by $(T_1,
Y_1), \cdots, (T_k,Y_k)$.

\begin{lem}
\label{lem-square} \be \label{l2cond} \E_{T}(\sup_{t\leq
T}|X_t|^2)<+\infty. \ee
\end{lem}

\begin{dem}
As usual,  we localize via  $\tau_n=\inf\{t>0, |X_t|\geq n\}$. Then,
for $t\leq T$,
\begin{align*}
    \E_{T}\left(\sup_{s\leq t}|X_{t\wedge \tau_n}|^2\right)&\leq C
\bigg(|X_0|^2+\int_0^{t\wedge
\tau_n}\int_{-1}^1 |y|^2 \bigg(\E_{T}( |X_s|^2)+\sigma^2(0,P_s)\bigg)\beta_1(y)dyds\\&\phantom{\leq C
\bigg(} +
\sum_{i=1}^k
|Y_i|^2\left(\E_{T}(|X_{T_i-}|^2)+\sigma^2(0,P_{T_i})\right){\bf 1}_{T_i\leq t\wedge \tau_n}\bigg)\\
&\leq C \bigg(|X_0|^2+\sup_{u\in[0,T]}\sigma^2(0,P_u)\bigg(1+\sum_{i=1}^k
|Y_i|^2\bigg)+\int_0^{t\wedge
\tau_n} \E_{T}( |X_s|^2)ds \\&\phantom{\leq K
\bigg(}+\sup_{i=1}^k |Y_i|^2\int_0^{t\wedge
\tau_n} \E_{T}( |X_{s-}|^2)\int_{|y|>1}N_2(dy,ds) \bigg).
 \end{align*}
 At this point  we apply
Gronwall's Lemma in its generalized form (with respect to the
measure $ds+\int_{|y|>1}N_2(dy,ds)$, see,  for example,  Ethier and Kurz
\cite{Kurz:86}, p. 498). The result then follows.
\end{dem}

Let us now explain our strategy to prove Theorem \ref{density}
before giving the technical details. We are going to prove that
there exists an a.s. positive random
 variable $DX_T$ such that $\E_T(DX_T)<+\infty$, and $ \exists C \;\forall \phi\in
C^\infty_c(\R),\;$
\begin{equation}
\label{ineq}|\E_T(\phi'(X_T)DX_T)|\leq C\|\phi\|_\infty,
\end{equation}
where $C^\infty_c(\R)$ denotes the space of infinitely
differentiable functions with compact support on the real line.
Indeed, in the special case $DX_T=1$ this inequality  implies
that the conditional law of $X_T$ given $(X_0,N_2^T)$ admits
a density (see,  for example,  Nualart
 \cite{Nualart:95},  p.79). If $DX_T\neq 1$, the law of $X_T$ under
 ${\mathbb Q}_T=\frac{DX_T.\P_T}{\E_T(DX_T)}$ admits a density. But
 since $DX_T > 0$,  a.s., then ${\mathbb Q}_T$ is
 equivalent to $\P_T$,  and the conditional law of $X_T$ given
 $(X_0,N_2^T)$ still admits a
density. Of course, this implies that the law $P_T$ of
 $X_T$ admits a density.

 We will prove inequality \eqref{ineq} employing the  stochastic calculus of
 variations. Consider perturbed paths of the process on the
 time interval $[0,T]$ and  introduce a parameter $\lambda\in [-1,1]$,  sufficiently
 close to $0$.   The perturbed Poisson measure
 $N_1^\lambda$ will satisfy $N_1^0=N_1$,  and be such that,  for a well
 chosen $\P_T$-martingale $(G^\lambda_t)_{t\leq T}$,  the law of its restriction to $[0,T]$ under $G^\lambda_T.\P_T$  will be
 equal
 to the law of $N_1$ under $\P_T$.
 The process $X^\lambda$ will be defined like $X$, only replacing $N_1$
 by $N_1^\lambda$ in the stochastic differential equation. Then, for sufficiently
 smooth functions $\phi$, we will have
 \be
 \label{loi}
 \E_T(\phi(X_T))=\E_T(G^\lambda_T\phi(X^\lambda_T)).
 \ee
 Differentiating in $\lambda$ at $\lambda=0$,  in a  sense that we yet have
 to define, we will obtain
\be
\label{intpartie}\E_T(\phi'(X_T)DX_T)=-\E_T(DG_T\phi(X_T)),
\ee
where $DX_T=\frac{d}{d\lambda}X^\lambda_T|_{\lambda=0}$,  and
$DG_T=\frac{d}{d\lambda}G^\lambda_T|_{\lambda=0}$. Then one easily
deduces
 \eqref{ineq} with $C=\E_T(|DG_T|)$.

 Let us describe the perturbation we are interested in.
Let $g$ be an increasing function  of class
$C^{\infty}_b(\mathbb{R})$, equal to $x$ on $[-{1\over 2}; {1\over
2}]$, equal to $1$, for $x\geq 1$, and to $-1$, for $x\leq -1$. Note
that $\|g\|_\infty\leq 1$, and that $g(x)x>0$,  for $x\in
\mathbb{R}^*$.

Now,  we define the predictable function $v: \Omega\times
[0,T]\times
  [-1,1] \mapsto \mathbb{R}$ via the formula
   \be
    \label{v}
v(s,y)={\bf 1}_{\{s>S\}}(1+y\sigma'_x(X_{s-},P_s))\ g(\sigma(X_{s-},P_s))\ k(y)\,
\ee
where $S$ is a stopping time that we are going to choose later on in order to ensure
that $DX_T>0$,  a.s. (see the discussion before Proposition \ref{positivity}).
It is easy to verify that the function
 $y\mapsto v(t,y)$ is of class $C^1$ on $[-1,1]$, and in what follows we will denote its derivative by   $v'(t,y)$, Also, for
  every $\omega,t,$ and $y$,
  \be
  \label{prop-v}
  |v(t,y)|\leq (1+K_1)k(y)\quad {\rm and } \quad |v'(t,y)|\leq
  K_1k(y)+(1+K_1)|k'(y)|.
  \ee


For  $\lambda\in [-1,1]$, let us introduce the perturbation
function
 \be
  \label{perturbation}
\gamma^\lambda(t,y)=y+\lambda v(t,y). \ee We can easily check
that,  for every $\omega,$ and $t$, the map $\, y\mapsto
\gamma^\lambda(t,y)$ is an increasing bijection from $[-1,1]$ into
itself, since by \eqref{hypbound} and \eqref{prop-v}, $|v'|\leq
{1\over 2}$, and $k(-1)=k(1)=0$.

Let us denote by $N_1^\lambda$ the image measure of the
Poisson point measure $N_1$ via the mapping $\gamma^\lambda$ defined,
 for any Borel subset $A$ of $[0,T]\times[-1,1]$, by the integral
$$
N_1^\lambda(A)=\int{\bf 1}_A(t,\gamma^\lambda(t,y))N_1(dy,dt).
$$
 We also introduce the  function
\be \label{densite-girsanov} V^\lambda(s,y)={\beta_1(y+\lambda
v(s,y))\over \beta_1(y)}(1+\lambda v'(s,y)), \ee which appears
below in the definition of the process $G^\lambda$ (in Proposition
\ref{def-g}). As a preliminary step, we obtain the following
estimates concerning $V^\lambda$.

\begin{lem}
\label{V} There exists a constant $C$ such that,  for almost all
$\omega$, and for all $s \in [0,T]$,  \begin{align} &\sup_{\lambda, y \in[-1,1]}{1\over
\lambda}|V^\lambda(s,y)-1| \leq {1\over 2}\
,\label{borneV}\\
&\sup_{\lambda\in[-1,1]}{1\over
\lambda^2}\int_{-1}^1|V^\lambda(s,y)-1|^2 \beta_1(y)dy
\leq C,\ \label{derpremV}\\
&\sup_{\lambda\in[-1,1]}{1\over
\lambda^4}\int_{-1}^1\left|V^\lambda(s,y)-1-\lambda{d\over
d\lambda}V^\lambda(s,y)/_{\lambda=0}\right|^2 \beta_1(y)dy \leq C. \label{dersecV}
\end{align}
\end{lem}

\begin{dem}
  Inequality \eqref{borneV} follows from \eqref{hypborne}. Also, one has
\be
 {d\over d\lambda}V^\lambda(s,y)&=& v'(s,y){\beta_1(y+\lambda v(s,y))
 \over \beta_1(y)}+{\beta_1'(y+\lambda v(s,y))\over
 \beta_1(y)}v(s,y)(1+\lambda v'(s,y)),\nonumber\\
 \ \label{derprem}
 \ee
 and
 \be
 {d^2\over d\lambda^2}V^\lambda(s,y)&=&v^2(s,y)(1+\lambda v'(s,y)){\beta''_1(y+\lambda v(s,y))
 \over \beta_1(y)}+2{\beta_1'(y+\lambda v(s,y))\over
 \beta_1(y)}v(s,y) v'(s,y).\nonumber\\
 \ \label{dersec}
 \ee
Since,  for $\lambda\in[-1,1]$, we have estimates,
$$
\left|\frac{V^\lambda(s,y)-1}{\lambda}\right|^2\leq
2\left|{d\over
d\lambda}V^\lambda(s,y)/_{\lambda=0}\right|^2+\frac{2}{\lambda^2}\left|V^\lambda(s,y)-1-\lambda{d\over
d\lambda}V^\lambda(s,y)/_{\lambda=0}\right|^2,
$$
and
$$
  \left|V^\lambda(s,y)-1-\lambda{d\over
d\lambda}V^\lambda(s,y)/_{\lambda=0}\right|^2\leq
\frac{\lambda^4}{4}\sup_{\mu\in[-1,1]}\left|{d^2\over
    d\mu^2}V^\mu(s,y)\right|^2,
$$
one deduces \eqref{dersecV} (resp. \eqref{derpremV}) from
\eqref{hypderive},  and \eqref{hypbound}(resp. \eqref{hypbeta}, \eqref{hypderive},  and \eqref{hypbound}).
\end{dem}

We are now ready to introduce the promised earlier definition of
the process  $G^\lambda$.

\begin{prop}
\label{def-g} (i) For every $\lambda\in [-1,1]$, the stochastic
differential equation \be \label{mart-girsanov} G^\lambda_t=
1+\int_{(0,t]\times \mathbb{R}} G^\lambda_{s-}
(V^\lambda(s,y)-1)\tilde{N}_1(dy,ds), \ee has a unique solution
$G^\lambda$ which is a strictly positive martingale under $\P_T$ and
such that
\be
 \label{G-L2}
\sup_{\lambda\in[-1,1]}\E_T\left(\sup_{t\leq
T}|G^\lambda_t|^2\right)<+\infty.
\ee
 (ii) The law of $\ N_1^\lambda$
under $\ \mathbb{P}_T^\lambda=G^\lambda_T.\mathbb{P}_T$ is the
same as the law of $N_1$ under $\mathbb{P}_T$.
\end{prop}

\begin{dem}
{\it (i)}  Thanks to \eqref{derpremV}, the
stochastic integral
$M^\lambda_t=\int_{(0,t]\times\R}(V^\lambda(s,y)-1)\tilde{N}_1(dy,ds)$
is well defined and is a $\P_T$ square integrable martingale. The
unique solution to \eqref{mart-girsanov} is the exponential
martingale $\ G^\lambda_t={\cal
E}(M^\lambda)_t=e^{M^\lambda_t}\prod_{0<s\leq t}(1+\Delta
M^\lambda_s)e^{-\Delta M^\lambda_s}\ $ given by the Dol\'eans-Dade
formula.

Using \eqref{borneV}, we remark that the jumps of $M^\lambda$
are more than $-1/2$ so that  $G^\lambda_t$ is positive for each
$t$. Moreover, using \eqref{derpremV}, Doob's inequality and
Gronwall's Lemma, we deduce from \eqref{mart-girsanov} that \eqref{G-L2} holds.

{\it (ii) } Let us denote $\mu(dy,dt)={\bf 1}_{|y|\leq
      1}\beta_1(y)dy dt$,  and compute the image measure $\gamma^\lambda(V^\lambda.\mu)$. For
a Borel subset $A$ of $[0,T]\times[-1,1]$, we have \be \gamma^\lambda(V^\lambda.\mu)(A)&=&
\int {\bf 1}_{A}(t,y+\lambda v(t,y))
V^\lambda(t,y)\beta_1(y)dy  dt\nonumber\\
&=& \int  {\bf 1}_{A}(t,y'){\beta_1(y')\over
\beta_1(y)}(1+\lambda v'(t,y))\beta_1(y){dy'\over 1+\lambda
v'(t,y)}
  dt\nonumber\\
&=&\int  {\bf 1}_{A}(t,y')\beta_1(y')dy' dt = \mu(A), \ee
where $y'=y+\lambda v(t,y)$. Hence \be \label{image}
\gamma^\lambda(V^\lambda.\mu)=\mu. \ee Since $N_1$ is independent
from $(X_0,N_2^T)$, the compensator of $N_1$ under $\P_T$ is
$\mu$. By the Girsanov's theorem for random measures (cf.
Jacod and Shiryaev \cite{Jacod:87},  p. 157), its compensator under $\
\mathbb{P}^\lambda_T=G^\lambda_T.\mathbb{P}_T$ is $\
V^\lambda.\mu\ $ and thus, the compensator  of $\
N_1^\lambda=\gamma^\lambda(N_1)\ $ is equal to $\
\gamma^\lambda(V^\lambda.\mu)=\mu$. We have thus proven that the
compensator of $\ N_1^\lambda$ under $\ \mathbb{P}^\lambda_T$ is
$\ \mu$,  and the second assertion in the proposition follows.
 \end{dem}

Next, we study the differentiability of $G^\lambda$ with respect to
the parameter $\lambda$,  at $\lambda=0$.
\begin{prop} (i)
 The process
  \be
 \label{dg}
 DG_t&=&\int_{(0,t]\times \R}{d\over
 d\lambda}V^\lambda(s,y)_{/\lambda=0}\tilde{N}_1(dy,ds)
 = \int_{(0,t]\times \mathbb{R}} \left(v'(s,y)+{\beta_1'(y)\over
 \beta_1(y)}v(s,y)\right)\tilde{N}_1(dy,ds)\nonumber\\
 && \
 \ee
  is well defined, and such that
  \be
   \label{dg-fini}\E_T(\sup_{t\leq T}|DG_t|^2)<+\infty.
    \ee

(ii) The process $G^\lambda$ is $L^2$-differentiable at $\lambda=0$,
 with the derivative $DG$ which is understood in the following sense:

 \be
 \label{der-g}
 \E_T\left(\sup_{t\leq
 T}|G^\lambda_t-1-\lambda DG_t|^2\right)=o(\lambda^2), \quad {\rm a.s.},
 \ee
 as  $\lambda$ tends to $0$.
\end{prop}

\begin{dem} {\it (i) }
 Thanks to \eqref{derprem} and \eqref{hypbeta}, for almost all $\omega$, and
 all $s\in[0,T]$,
\begin{equation}
   \int_{-1}^1\left|{d\over
       d\lambda}V^\lambda(s,y)_{/\lambda=0}\right|^2\beta_1(y)dy\leq C\label{majodv}.
\end{equation}
Therefore
the process $DG_t$ is well
defined and satisfies \eqref{dg-fini}.

{\it (ii)}  Moreover, one has
\ben &&\E_T\left(\sup_{t\leq
 T}|G^\lambda_t-1-\lambda DG_t|^2\right)  \\
& \leq& C \int_{(0,t]\times
\mathbb{R}}\E_T\left(\left( G^\lambda_{s}
\left(V^\lambda(s,y)-1\right)-\lambda {d\over d\lambda}V^\lambda(s,y)_{/\lambda=0}\right)^2\right)\beta_1(y)dyds\\
&\leq &C \int_{(0,t]\times \mathbb{R}}\E_T\left(\left(
G^\lambda_{s}
\left(V^\lambda(s,y)-1-\lambda {d\over d\lambda}V^\lambda(s,y)_{/\lambda=0}\right)\right)^2\right)\beta_1(y)dyds\\
&&\hskip 1cm+\lambda^2 C \int_{(0,t]\times
\mathbb{R}}\E_T\left(\left( {d\over
d\lambda}V^\lambda(s,y)_{/\lambda=0})\right)^2(G^\lambda_s-1)^2\right)\beta_1(y)dyds.
\een
Now,  according to   \eqref{dersecV} and \eqref{G-L2}, we
obtain that \ben &&\int_{(0,t]\times \mathbb{R}}\E_T\left(\left(
G^\lambda_{s} (V^\lambda(s,y)-1-\lambda {d\over
d\lambda}V^\lambda(s,y)_{/\lambda=0})\right)^2\right)\beta_1(y)dyds\leq C \lambda^4 t.
 \een
 Furthermore, by \eqref{majodv},
$$
\int_{(0,t]\times \mathbb{R}}\E_T\left(\left(
{d\over
d\lambda}V^\lambda(s,y)_{/\lambda=0})\right)^2(G^\lambda_s-1)^2\right)\beta_1(y)dyds
\leq \int_0^t \E_T((G^\lambda_s -1)^2)ds,
$$
 and
using \eqref{derpremV} and \eqref{G-L2}, we may show that, for each
$t\leq T$,
$$\E_T((G^\lambda_t
-1)^2)=\int_{(0,t]\times\R}\E_T((G^\lambda_s(V^{\lambda}(s,y)-1))^2)\beta_1(y)dyds
\leq C \lambda^2.$$ This  concludes the proof.
\end{dem}

In the next step we  define the perturbed stochastic differential equation.
Let us recall that the probability measures $P_t$ are fixed and
are considered as time-dependent  parameters. Thus the process $X$
is a function $F_P(X_0,N_1, N_2)$ of the triplet $(X_0,N_1, N_2)$.

Define $X^\lambda: =F_P(X_0,N_1^\lambda,
N_2)$. Hence,  using Proposition \ref{def-g} {\it (ii)}, the law of
$X^\lambda$ under $\mathbb{P}^\lambda_T$ is equal to the law of
$X$ under $\mathbb{P}_T$. A simple computation shows that
$X^\lambda$ is a
 solution of the stochastic differential equation
  \be
\label{eds-perturb} X^\lambda_t&=&X_0+\int_{(0,t]\times
\mathbb{R}} y\ \sigma(X^\lambda_{s-},P_s)\
\left(N_1^\lambda(dy,ds)-\beta_1(dy)ds\right) +\int_{(0,t]\times
\mathbb{R}}
y\ \sigma(X^\lambda_{s-},P_s)N_2(dy,ds)\nonumber\\
&=&X_0+ \int_{(0,t]\times \mathbb{R}} (y+\lambda
v(s,y))\ \sigma(X^\lambda_{s-},P_s)\ \left(N_1(dy,ds)-V^\lambda(s,y)\beta_1(dy)ds\right)\nonumber\\
&&\hskip 3cm  +\int_{(0,t]\times \mathbb{R}} y\
\sigma(X^\lambda_{s-},P_s)N_2(dy,ds),
\quad \quad \quad\qquad  ( \hbox{ since } \gamma^\lambda(V^\lambda.\mu)=\mu) \nonumber\\
&=&X_0+\int_{(0,t]\times \mathbb{R}}
y\sigma(X^\lambda_{s-},P_s)\tilde{N}_1(dy,ds)+\lambda\int_{(0,t]\times
\mathbb{R}}
\sigma(X^\lambda_{s-},P_s)v(s,y)\tilde{N}_1(dy,ds)\nonumber\\&+&\int_{(0,t]\times
\mathbb{R}} \hskip -5mm y\
\sigma(X^\lambda_{s-},P_s)N_2(dy,ds)-\int_{(0,t]\times
\mathbb{R}}\hskip-5mm (y+\lambda v(s,y))
\sigma(X^\lambda_{s-},P_s)(V^\lambda(s,y)-1)\beta_1(y)dyds.\nonumber\\
\ee
 Using \eqref{hypint} for the second term,  the fact that
$\int_{-1}^1(y^2+k^2(y))\beta_1(y)dy<+\infty$, \eqref{derpremV},
and the Cauchy-Schwarz inequality for the last term, we easily prove
that equation \eqref{eds-perturb} has a unique
 pathwise solution.

Let us now show  that $X^\lambda$ is differentiable in $\lambda$,
at $\lambda=0$, in the  $L^2$-sense. More precisely we  have the following
\begin{prop}
\be
 (i)&&\label{cont-x} \E_T\left(\sup_{t\leq
T}|X^\lambda_t - X_t|^4\right)\leq C \lambda^4.\\
\label{diff-x} (ii)&& \E_T\left(\sup_{t\leq T}|X^\lambda_t -
X_t-\lambda DX_t|^2\right)= o(\lambda^2),
\ee
as  $\lambda$ tends
to $0$, where $DX$ is a solution of the affine stochastic differential
equation
$$
  \label{dx} DX_t= \int_{(0,t]\times \mathbb{R}}
y\sigma'_x(X_{s-},P_s) DX_{s-}\ \tilde{N}_1(dy,ds)
+\int_{(0,t]\times \mathbb{R}}
\sigma(X_{s-},P_s)v(s,y)\tilde{N}_1(dy,ds) \nonumber
$$
\be
 +\int_{(0,t]\times \mathbb{R}}\hskip -7mm y\sigma'_x(X_{s-},P_s) DX_{s-}\
N_2(dy,ds) -\int_{(0,t]\times [-1,1]} \hskip -10mm y\
\sigma(X_{s-},P_s)(\beta_1(y)v'(s,y)+\beta'_1(y)v(s,y))dyds.\nonumber\\
\ \ee
\end{prop}

\begin{dem}
 In order to prove assertion {\it  (i)}, we need the
following moment estimate  :
\begin{equation}\sup_{\lambda\in[-1,1]}\E_T\left(\sup_{t\leq
    T}|X^\lambda_t|^4\right)<+\infty\label{mom4per}.\end{equation}
It relies on the following classical estimation (see,
\cite{Bichteler:83},  Lemme (A.14)):
\begin{align}
   \E_T\left(\left(\int_{(0,t]\times \mathbb{R}}
H_s \rho(y)\tilde{N}_1(dy,ds)\right)^4\right)\leq
C&\left(\left(\int_{-1}^1\rho^2(y)\beta_1(y)dy\right)^2+\int_{-1}^1\rho^4(y)\beta_1(y)dy\right)\notag\\&\times\int_0^t
\E_T\left(\sup_{u\leq s}|H_u|^4\right)ds,\label{mom4}
\end{align}
  for any
predictable process $H$,  and any measurable function
$\rho:[-1,1]\mapsto \R$ such that the right-hand side is
finite. Conditioning by $\ N_2^T$, the times and the amplitudes of  jumps of
$N_2$ on $(0,t]$ are given by $(T_1, Y_1), \cdots, (T_k,Y_k)$,  and
\begin{align*}
   \E_{T}\bigg(\bigg|\int_0^ty\sigma(X^\lambda_{s-},P_s)&N_2(dy,ds)\bigg|^4\bigg)
\leq
 C\sum_{i=1}^k
Y_i^4\left(\E_{T}(|X^\lambda_{T_i-}|^4)+\sigma^4(0,P_{T_i})\right)\\
&\leq C\sup_{i=1}^k
|Y_i|^4\left(\sup_{u\in[0,T]}\sigma^4(0,P_{u})+\int_0^t\int_{|y|>1}\E_{T}(|X^\lambda_{s-}|^4)
N_2(dy,ds)\right).
\end{align*}
Applying \eqref{mom4} with $\rho(y)=y$,  and
$\rho(y)=k(y)$, \eqref{derpremV} and Gronwall's Lemma with respect
to the measure $ds+\int_{|y|>1}N_2(dy,ds)$, we easily check
\eqref{mom4per} and deduce
$$\sup_{\lambda\in[-1,1]}\E_T\left(\sup_{t\leq
    T}|\sigma(X^\lambda_t,P_t)|^4\right)<+\infty.$$
Now,  we write $X^\lambda_t-X_t$ using \eqref{eds-stable} and
\eqref{eds-perturb}. Assertion {\it (i)} is obtained  following an analogous argument.


To prove {\it (ii) } we need to isolate the term $Z^\lambda_t=X^\lambda_t - X_t-\lambda
DX_t$, and as in \cite{Bichteler:83}, Theorem (A.10), we write
\ben
&&X^\lambda_t -
X_t-\lambda DX_t\\
&=&\int_{(0,t]\times \mathbb{R}}y\ Z^\lambda_{s-}
\sigma'_x(X_{s-},P_s)(\tilde{N}_1(dy,ds)+N_2(dy,ds))\\
&&+ \int_{(0,t]\times \mathbb{R}} \hskip -5mm y\
\left(\sigma(X^\lambda_{s-},P_s)-\sigma(X_{s-},P_s)-\sigma'_x(X_{s-},P_s)(X^\lambda_{s-}-X_{s-})\right)
(\tilde{N}_1(dy,ds)+N_2(dy,ds))\\
&& +\int_{(0,t]\times
\mathbb{R}}\lambda v(s,y)\left(\sigma(X^\lambda_{s-},P_s)-\sigma(X_{s-},P_s)\right)\tilde{N}_1(dy,ds)\\
&&-\int_{(0,t]\times \mathbb{R}}y\ \sigma(X_{s},P_s)
\left(V^\lambda(s,y)-1-\lambda {d\over
d\lambda}V^\lambda(s,y)_{/\lambda=0}\right)\beta_1(y)dy ds\\
&&-\int_{(0,t]\times \mathbb{R}}y\
\left(\sigma(X^\lambda_{s},P_s)-\sigma(X_{s},P_s)\right)(V^\lambda(s,y)-1)\beta_1(y)dy ds\\
&&-\int_{(0,t]\times \mathbb{R}}\lambda\ v(s,y)
\sigma(X^\lambda_{s},P_s)(V^\lambda(s,y)-1)\beta_1(y)dy ds.
 \een

 Under Hypotheses ({\bf H}) and ({\bf H}$_1$) all the integral terms, except
 the first one, are of  order $\lambda^2$.
 Indeed, for the second term, we  use   Taylor's
  expansion  of $\sigma$, and \eqref{cont-x}; for the third term, we use \eqref{hypint},
   and \eqref{cont-x}; for the fourth
  we use the Cauchy-Schwarz inequality, and \eqref{dersecV}; for the fifth term, we use the
   Cauchy-Schwarz inequality, \eqref{derpremV}, and \eqref{cont-x}; for
  the sixth term, we use Cauchy-Schwarz inequality,
  \eqref{derpremV},  and \eqref{hypint}.
  Then, as previously, using  Gronwall's Lemma for the conditional expectation, we obtain the result.
\end{dem}

The term $DX_t$ requires our special attention. Observe that,  after
integration by parts (in the variable $y$),  the last term in
\eqref{dx} cancels the compensated part of
$$
\int_{(0,t]\times
\mathbb{R}} \sigma(X_{s-},P_s) v(s,y) \tilde{N}_1(dy,ds),
$$
 and one
obtains \be \label{dv} DX_t=\int_0^t DX_{s-}dK_s + L_t \ee where
\be K_t&=&\int_{(0,t]\times \mathbb{R}}y\
\sigma'_x(X_{s-},P_s)\tilde{N}_1(dy,ds)+\int_{(0,t]\times
\mathbb{R}}y\
\sigma'_x(X_{s-},P_s) N_2(dy,ds), \label{kt}\\
L_t&=&\int_{(0,t]\times \mathbb{R}} \sigma(X_{s-},P_s) \ v(s,y)\
N_1(dy,ds).\label{lt} \ee As in \cite{Jacod:87},  Theorem 4.61, p.
59,  or in \cite{Bichteler:83}, we can solve  \eqref{dv}
explicitely. The jumps $\Delta K_s$ are of the form $\
y\sigma'_x(X_s,P_s)$. Thus $1+\Delta K_s$ may be equal to $0$ and
then the Doleans-Dade exponential
$${\cal E}(K)_t=e^{K_t}\Pi_{0<s\leq t}(1+\Delta K_s)e^{-\Delta K_s}$$ vanishes
from the first time when $\Delta K_s=-1$. We follow
\cite{Bichteler:83} to show that $DX_T\neq 0$, but  the strict
positivity (which has not been proved in the latter) necessitates
a careful analysis.

Let us define the sequence of stopping  times $S_1=\inf\{t>0,
\Delta K_t\leq -1\}$, $S_k=\inf\{t>S_{k-1}, \Delta K_t\leq -1\}$,
$S_0=0$. Since $\sigma'_x$ is bounded, there
is a finite number of big jumps on the time interval $[0,T]$,  so that
there exists an $n$ such that $S_n\leq T <S_{n+1}=\infty$,  and
$\P(S_n=T)=0$.

Solving equation  \eqref{dv}  gives
 \be  DX_t=&&{\cal
E}(K-K^{S_k})_t \bigg(DX_{S_k} +\int_{(S_k,t]}(1+\Delta
K_s)^{-1}{\cal E}(K-K^{S_k})_{s-}^{-1} dL_s\bigg)\nonumber\\&&\hbox{ if } S_k\leq t<S_{k+1},  \hbox{ and } t\leq T, \ee where
$K^{S_k}_t=K_{S_k\wedge t}$. In particular,
$$
 DX_T={\cal
E}(K-K^{S_n})_T \bigg(DX_{S_n} +\int_{(S_n,T]}(1+\Delta
K_s)^{-1}{\cal E}(K-K^{S_n})_{s-}^{-1} dL_s\bigg).
$$
Because of the definition of $S_n$, the exponential martingale ${\cal
  E}(K-K^{S_n})_{s}$ is non-negative on $[S_n,T]$. If the perturbation
$\ v$ did not vanish before time $S_n$, it would not be clear how
to control the sign of $DX_{S_n}$. That is why we choose $S=S_n$
in \eqref{v} :
$$v(s,y)={\bf 1}_{\{s>S_n\}}(1+y\ \sigma'_x(X_{s-},P_s))\ k(y)\
g(\sigma(X_{s-},P_s))$$
 so that $DX_{S_n}=0$. For this choice, we
obtain

\begin{prop}
\label{positivity} We have $DX_T>0$,  almost surely.
\end{prop}

\begin{dem}
One has $DX_T={\cal
E}(K-K^{S_n})_TY_n$,  where $Y_n=\int_{(S_n,T]}(1+\Delta
K_s)^{-1}{\cal E}(K-K^{S_n})_{s-}^{-1} dL_s$.

Since $N_1$ and $N_2$ are independent, the sets of jumps are
almost surely distinct and then $1+\Delta K_s$ can be replaced by
$1$ every time the jump of $K$ comes from a jump of $N_2$.
 So,
 $$
 Y_n=\int_{(S_n,T]\times [-1,1]}h_n(s,y)N_1(ds,dy),
 $$
 where
 \be h_n(s,y)&=&{\cal
E}(K-K^{S_n})_{s-}^{-1}
(1+y\sigma'_x(X_{s-},P_s))^{-1}v(s,y)\sigma(X_{s-},P_s)\nonumber\\&=&{\cal
E}(K-K^{S_n})_{s-}^{-1}k(y)\sigma(X_{s-},P_s)g(\sigma(X_{s-},P_s))\geq
0.\nonumber\ee Let us consider the set $A_n=\{(\omega,s,y),
h_n(\omega,s,y)>0\}$ and define the stopping time
$\tau=\inf\{t>S_n, \int_{(S_n,t]\times [-1,1]}h_n(s,y) N_1(ds,dy)
>0\}=\inf\{t>S_n, \int_{(S_n,t]\times [-1,1]}{\bf
1}_{A_n}(s,y)N_1(ds,dy)>0\}.$

Using the definitions of $v$ and $S_{n+1}$, one knows that,  if
$S_n(\omega)<s\leq T\wedge S_{n+1}(\omega)$, then
$$
(\omega,s,y)\in A_n \Leftrightarrow \sigma(X_{s-}(\omega),P_s)
\neq 0,
$$
which is always the case in view of  Hypothesis ({\bf H}).

On the other hand, $\int_{(S_n,\tau]\times [-1,1]}{\bf
1}_{A_n}(s,y)N_1(ds,dy)~\leq~1$,\\
 so
$\ \E\left(\int_{(S_n,\tau]\times [-1,1]}{\bf
1}_{A_n}(s,y)N_1(ds,dy)\right)~\leq~1$  and,  for $\omega$ in a set
of probability $1$,
$$\int_{(S_n,\tau]\times [-1,1]}{\bf
1}_{A_n}(s,y) \beta_1(y)dy ds<+\infty.$$ These two remarks, and
the fact the $\ \int_{-1}^1 \beta_1(y)dy = +\infty$,  imply that
$\tau=S_n$,  almost surely. So  $Y_n$ is strictly
positive. Therefore $DX_T$ is strictly positive as well and the proof is complete.
\end{dem}

We are now in a position  to complete the proof of Theorem \ref{density}.

 \noindent {\bf Proof of Theorem 10\;:}  For $\phi\in
C^\infty_b(\mathbb{R})$, we differentiate the
expression \eqref{loi} in the $L^2$-sense  with respect to $\lambda$, at  $\lambda=0$,
and hence, we obtain the "integration-by-parts" formula
\eqref{intpartie}.  Then, since $\E_T|DG_T|<+\infty$, we obtain
\eqref{ineq}, which concludes the proof.
\hfill\mbox{\rule{2 true mm}{3 true mm}}

\section{The case of a symmetric stable driving process $Z$}
In this section, we assume that the real-valued driving process
$Z$ is a symmetric stable process with index $\alpha\in (0,2)$,
i.e., $Z$ is given by \eqref{Levy-dens} with
$\beta(dy)=\frac{K}{|y|^{1+\alpha}}dy$,  where $K>0$ is a
normalization constant. The generator of this process is the
fractional Laplacian (or, fractional symmetric derivative) of
order $\alpha$ on $\R$ :
$$D_{x}^\alpha f(x)=K\int_{\R} \left (f(x+y)-f(x)-{\bf 1}_{\{|y|\leq 1\}}f'(x)y
\right )\frac{dy}{|y|^{1+\alpha}}.$$
This operator may be defined alternatively via the Fourier transform
${\cal F}$ :
\ben
D_{x}^\alpha v(x)= K'{\mathcal F}^{-1}\Bigl (|\xi|^\alpha{\mathcal
  F}(v)(\xi)\Bigr )(x),\mbox{ with }K'>0.\een
When $\sigma:\R\times{\mathcal
  P}(\R)\rightarrow \R$ satisfies Hypotheses ({\bf H}), it is possible to explicitly calculate  the adjoint
${\mathcal L}^*[\nu]$  involved in the nonlinear Fokker-Planck equation
\eqref{fokplan}. For smooth functions
$\varphi,\psi:\R\mapsto\R$, one has
$$
\int_{\R}{\mathcal L}[\nu]\varphi(x)\psi(x)dx=K\int_{\R^2}
\left (\varphi(x+s(x)y)-\varphi(x)-{\bf 1}_{\{|y|\leq 1\}}s(x)y\varphi'(x)
\right)\frac{dy}{|y|^{1+\alpha}}\psi(x)dx,
$$
where $s(x)=\sigma(x,\nu)$. Setting $z=-s(x)y$,
and observing  that $\int_{\R}\left({\bf 1}_{\{|z|\leq s(x)\}}-{\bf
1}_{\{|z|\leq
    1\}}\right)\frac{zdz}{|z|^{1+\alpha}}=0$, one gets  \begin{align*}
      \int_{\R}\bigg(\varphi(x+s(x)y)-\varphi(x)&-{\bf 1}_{\{|y|\leq 1\}}s(x)y\varphi'(x)
\bigg)\frac{dy}{|y|^{1+\alpha}}\\&=\int_{\R}\left
(\varphi(x-z)-\varphi(x)+{\bf 1}_{\{|z|\leq 1\}}z\varphi'(x)
\right)|s(x)|^\alpha\frac{dz}{|z|^{1+\alpha}}.
    \end{align*}
Since
$$
\int_\R \varphi(x-z)[|s|^\alpha\psi](x)dx=\int_\R
\varphi(x)[|s|^\alpha\psi](x+z)dx,
$$
and
 $$
 \int_\R
\varphi'(x)[|s|^\alpha\psi](x)dx=-\int_\R
\varphi(x)[|s|^\alpha\psi]'(x)dx,
$$
invoking  Fubini's
theorem one concludes that
\begin{align*}
  \int_{\R}{\mathcal L}[\nu]\varphi(x)\psi(x)dx=K\int_{\R}\varphi(x)\int_{\R} \bigg([|s|^\alpha\psi](x+z)-[|s|^\alpha\psi](x)-{\bf 1}_{\{|z|\leq
    1\}}z[|s|^\alpha\psi]'(x)\bigg)\frac{dz}{|z|^{1+\alpha}}dx.
\end{align*}
Therefore
$$
{\mathcal L}^*[\nu]\psi(x)=D_{x}^\alpha(|\sigma(.,\nu)|^\alpha\psi(.))(x).
$$
Moreover,   the absolute continuity result given in Theorem
\ref{density} permits us to prove  existence of a function solution to the
nonlinear Fokker-Planck equation.

\begin{thm}
Let $m\in{\mathcal P}(\R)$,  and $\alpha\in(0,2)$. Assume that the function $\sigma(x,\nu)$ satisfies
   hypotheses {\bf (H)} and is Lipschitz continuous in its second variable when
   ${\mathcal P}(\R)$ is endowed with the modified Vaserstein metric $d_1$.
Then,  there exists
   a function $(t,x)\in(0,T]\times\R\mapsto p_t(x)\in\R_+$ such that,  for
   each $t\in(0,T]$, $\int_\R p_t(x)dx=1$ and, in the weak sense,
\begin{equation}\label{pdenonlin}
   \begin{cases}
      \partial_tp_t(x)=D_{x}^\alpha
(|\sigma(.,p_t)|^\alpha p_t(.))(x)\\
\lim_{t\rightarrow 0^+}p_t(x)dx=m(dx)\mbox{ for the weak convergence},
   \end{cases}
\end{equation}
where, by a slight abuse of notation, $\sigma(.,p_t)$ stands for   $\sigma(.,p_t(y)dy)$.
\end{thm}

\begin{dem}
Existence of a measure solution $(P_t)_{t\in[0,T]}$ to the nonlinear
Fokker-Planck equation follows from Propositions \ref{proppdmgene} and
\ref{propedsgene}. So to conclude the proof, it is enough to exhibit a perturbation function
$k(y)$ satisfying hypotheses ({\bf H}$_1$) with
$\beta_1(y)={\bf 1}_{\{|y|\leq 1\}}\frac{K}{|y|^{1+\alpha}}$.
Then, by Theorem \ref{density},  for each $t\in(0,T]$, we have $P_t=p_t(x)dx$ .

For $\gamma>\frac{\alpha}{2}$,  and $\varepsilon\in(0,1/2)$, let $k_\varepsilon$ denote the even
function on $[-1,1]$ defined by
\begin{equation*}
   k_\varepsilon(y)=
   \begin{cases}
      y^{1+\gamma}, &\mbox{ for } y\in[0,\varepsilon],\\
\varepsilon^{1+\gamma}+(1+\gamma)\varepsilon^\gamma(y-\varepsilon)-(1+c)(y-\varepsilon)^{1+\gamma},&\mbox{
  for } y\in[\varepsilon,2\varepsilon]  \\
     (1+\gamma-c)\varepsilon^{1+\gamma}-c(1+\gamma)\varepsilon^\gamma(y-2\varepsilon),
& \mbox{ for } y\in[2\varepsilon,1],
   \end{cases},
\end{equation*}
where
$c=\frac{(1+\gamma)\varepsilon}{(1+\gamma)-\varepsilon(1+2\gamma)}$,  so
that $k_\varepsilon(1)=0$. The function $k_\varepsilon$ is non-negative
and $C^1$ on
$[0,1]$, satisfies \eqref{hypint}  and,  $\forall y\in[-1,1],\;$
\begin{equation}
k_\varepsilon(y)\leq
(2+\gamma)\varepsilon^{1+\gamma},\quad |k'_\varepsilon(y)|\leq
(1+\gamma)\max(1,c)\varepsilon^{\gamma},\;\quad {\rm and } \quad \frac{k_\varepsilon(y)}{|y|}\leq
(1+\gamma)\varepsilon^\gamma.\label{majosk}\end{equation}
In particular,  for
 small enough $\varepsilon$, $k_\varepsilon$ satisfes \eqref{hypborne}. Since
\begin{align*}
   \left|\frac{\beta_1'(y)}{\beta_1(y)}(1+ay)k_\varepsilon(y)+ak_\varepsilon(y)+(1+a y)k_\varepsilon'(y)\right|^2\beta_1(y)
&\leq
C\left[\frac{k_\varepsilon^2(y)}{y^2}+k_\varepsilon^2(y)+(k_\varepsilon')^2(y)\right]\beta_1(y)\\
&\sim
C'|y|^{-(1+\alpha-2\gamma)},
\end{align*}
  in the neighbourhood of $0$,
\eqref{hypbeta} is satisfied as well. In the same way, in the neighbourhood of $0$,
\begin{align*}
  \sup_{a\in[-K_1,
K_1],\lambda\in[-1,1]}\bigg(&\left|{\beta_1''(y+\lambda(1+a y)k_\varepsilon(y))\over \beta_1(y)}
(1+\lambda(a k_\varepsilon(y)+(1+ay)k_\varepsilon'(y))\right|^2k_\varepsilon^2(y)\\&+\left|{\beta_1'(y+\lambda(1+a
y)k_\varepsilon(y))\over \beta_1(y)}\right|^2 \bigg)k_\varepsilon^2(y)\beta_1(y)
\leq
C\left(\frac{|y|^{2+2\gamma}}{y^{4}}+\frac{1}{y^2}\right)|y|^{1+2\gamma-\alpha},
\end{align*}
and \eqref{hypderive} is satisfied.

Finally, for $a\in[-K_1,
K_1]$ and $y,\lambda\in[-1,1]$, by \eqref{majosk}, for $\varepsilon<((1+K_1)(1+\gamma))^{-1/\gamma}$,
\begin{align*}
&\frac{1}{|\lambda|}\left|{\beta_1(y+\lambda(1+a
y)k_\varepsilon(y))\over \beta_1(y)}(1+\lambda(a
k_\varepsilon(y)+(1+ay)k_\varepsilon'(y))-1\right|\\
=&\frac{1}{|\lambda|}\frac{|1-|1+\lambda(1+ay)\frac{k_\varepsilon(y)}{y}|^{1+\alpha}+\lambda(a
k_\varepsilon(y)+(1+ay)k_\varepsilon'(y))|}{|1+\lambda(1+ay)\frac{k_\varepsilon(y)}{y}|^{1+\alpha}}\\
\leq&
\frac{1}{(1-(1+K_1)(1+\gamma)\varepsilon^\gamma)^{1+\alpha}}\Bigl[ (1+\alpha)(1+(1+K_1)(1+\gamma)\varepsilon^{\gamma})^{\alpha}(1+K_1)(1+\gamma)\varepsilon^{\gamma}\\
&\hskip 70mm +K_1(2+\gamma)\varepsilon^{1+\gamma}+(1+K_1)(1+\gamma)\max(1,c)\varepsilon^{\gamma}\Bigr],
\end{align*}
and \eqref{hypbound} is also satisfied for small enough $\varepsilon$.
\end{dem}

 \begin{rem}\rm
One of the motivations for our work was to generalize the probabilistic
approximation  of the porous medium equation
\begin{equation}
\partial_tp_t(x)=D^2_x( p_t^q (x)), \qquad q>1,
\label{porousmedium}
\end{equation}developed, among others, by Jourdain \cite{Jourdain:00} to
the fractional case where $D^2_x$ is replaced by $D^\alpha_x$.
The equation \eqref{porousmedium}, which describes percolation of gases through porous media, and which is usually derived by combining the power type equation of state relating pressure to gas density $p$, conservation of mass law, and so called Darcy's law describing the local gas velocity as the gradient of pressure, goes back, at least, to the 1930's (see, e.g., Muskat \cite{Muskat:37}). The major steps in the development of the mathematical theory of \eqref{porousmedium} were the discovery of the family of its self-similar solutions by Barenblatt (see, \cite{Barenblatt:52}, and \cite{Barenblatt:96}) who obtained this equation in the context of heat propagation at the initial stages of a nuclear explosion, and an elegant uniqueness result for \eqref{porousmedium} proved by Br\'ezis and Crandall \cite{BR:79}. A summary of some of the newer developments in the area of the standard porous medium equation can be found in a survey by Otto \cite{Otto:01}.

However, in a number of recent physical papers, an argument was made that some of the fractional scaling observed in flows-in-porous-media phenomena cannot be modeled in the framework of \eqref{porousmedium}. In particular, Meerschaert, Benson and Baeumer \cite{MBB:99} replace the Laplacian $D_{x}^2$ in \eqref{porousmedium} by the fractional Laplacian $D_{x}^\alpha$ while considering the linear case ($q=1$) in a multidimensional case of anomalous (mostly geophysical) diffusion in porous media, while Park, Kleinfelter and Cushman \cite{PKC:05} continue in this tradition and derive scaling laws and (linear) Fokker-Planck equations for 3-dimensional porous media with fractal mesoscale.

On the other hand, Tsallis and Bukman \cite{TB:96} suggest an alternative approach to the anomalous scaling problem (in porous media, surface growth, and certain biological phenomena) and consider an equation of the general form
\begin{equation}
\partial_tp_t^r(x)=-D_x(F(x)p_t^r(x)) +D^2_x( p_t^q (x)), \qquad r,q\in \R,
\label{TBequation}
\end{equation}
where $F(x)$ is an external force. The authors manage to find exact
solutions for this class of equations using ingeniously the concept of
Renyi (-Tsallis) entropy but, significantly, suggest in the conclusion
of their paper that it would be desirable to develop physically
significant models for which further unification can possibly be
achieved by considering the generic case of a {\it nonlinear}
Fokker-Planck-like equation with {\it fractional} derivatives. This is
what we endeavored to do taking as our criterion of "physicality" the
existence of an approximating interacting particle scheme. For the most
obvious, simply-minded generalization,
$\partial_tp_t(x)=D^\alpha_x(p_t^q(x))$, that physical interpretation
seems to be missing, or, at least, we were unable to produce it and, as
a result, our study lead us to settle on an equation like
\eqref{pdenonlin}. Indeed for
$$\sigma(x,\nu)=\left(g_\varepsilon*\nu(x)\right)^s\mbox{
  with
}\varepsilon>0,\;g_\varepsilon(x)=\frac{1}{\sqrt{2\pi\varepsilon}}e^{-\frac{x^2}{2\varepsilon}}\;
\mbox{and}\;s>0,$$ \eqref{pdenonlin} writes $\partial_t
p_t=D^\alpha_x((g_\varepsilon*p_t)^{\alpha s}p_t)$ which, for now,
we are viewing as a "physically justifiable", fractional, {\it
  and} strongly nonlinear ``extension" of the classical porous medium
equation. Of course,
 this is only the beginning of the effort to understand these types of models.

 \end{rem}

\end{document}